\documentclass[journal]{article}
\usepackage{titlesec}
\usepackage{amsmath}
\usepackage{amssymb}
\usepackage{multirow}
\usepackage{amsfonts}
\usepackage{mathrsfs}
\usepackage{booktabs}
\usepackage{authblk}
\usepackage{algorithmicx}
\usepackage{algorithm}
\usepackage{booktabs}
\usepackage{float}
\usepackage{algpseudocode}
\usepackage{amsthm}
\usepackage{url}
\usepackage{makecell}
\usepackage{latexsym}
\usepackage{graphicx}
\usepackage[table]{xcolor}
\usepackage[most]{tcolorbox}
\usepackage{rotating}
\usepackage{subfig}
\usepackage{color,varioref}
\usepackage{longtable}
\allowdisplaybreaks[4]
\usepackage[colorlinks,linkcolor=blue]{hyperref}
\def\keywords{\vspace{.5em}{\bf{Keywords}:\,\relax}}

\newtheorem{remark}{Remark}
\newtheorem{corollary}{Corollary}
\usepackage{geometry}
\begin{document}

\title{Efficient QR-Based CP Decomposition Acceleration via Restructured Dimension Tree  and  Customized Extrapolation}

\author{Wenchao Xie\thanks{School of Mathematics and Computational Science, Xiangtan University, Xiangtan, 411105, China. Email: youmengx@hotmail.com}, \hskip 0.2cm
Jiawei Xu\thanks{School of Mathematics and Computational Science, Xiangtan University, Xiangtan, 411105, China. Email: xujiawei@smail.xtu.edu.cn}, \hskip 0.2cm 
Zheng Peng\thanks{School of Mathematics and Computational Science, Xiangtan University, Xiangtan, 411105, China. Email: pzheng@xtu.edu.cn}, \hskip 0.2cm
Qingsong Wang\thanks{Corresponding author. School of Mathematics and Computational Science, Xiangtan University, Xiangtan, 411105, China. Email: nothing2wang@hotmail.com}
}
\date{}	
\maketitle
\begin{abstract}
The canonical polyadic (CP) decomposition is one of the most widely used tensor decomposition techniques. The conventional CP decomposition algorithm combines alternating least squares (ALS) with the normal equation. However, the normal equation is susceptible to numerical ill-conditioning, which can adversely affect the decomposition results. To mitigate this issue, ALS combined with QR decomposition has been proposed as a more numerically stable alternative. Although this method enhances stability, its iterative process involves tensor-times-matrix (TTM) operations, which typically result in higher computational costs. To reduce this cost, we propose restructured dimension tree, which increases the reuse of intermediate tensors and reduces the number of TTM operations. Compared with the standard dimension tree structure, this dimension tree structure can reduce the computational complexity of TTM operations for tensors of any order by 33\%. Additionally, we introduce a customized extrapolation strategy in the CP-ALS-QR algorithm, leveraging the unique structure of the matrix $\mathbf{Q}_0$ to further accelerate convergence. By integrating these two techniques, we propose a novel CP decomposition algorithm that significantly improves iteration efficiency, achieving up to twofold acceleration on datasets with certain specific structures. Numerical experiments on five real-world datasets show that, compared with the baseline algorithm, our proposed algorithm improves iteration efficiency while simultaneously enhancing fitting accuracy.
\end{abstract}

\keywords{CP decomposition, ALS, QR decomposition, Dimension tree, Extrapolation}

\maketitle

\section{Introduction}

Tensors, also known as multidimensional arrays, have been widely studied due to their ability to extend the concept of matrices and efficiently model multidimensional data.
This ability makes them a powerful tool in various applications, including signal processing \cite{Lathauwer1998,Nion2010}, computer vision \cite{A8}, recommendation systems \cite{Rendle2010,Symeonidis2008}, image processing \cite{Sonka2014}, knowledge bases \cite{Carlson2010}, and other domains \cite{Ji2016,TG}. 
 In these applications, tensor methods primarily rely on tensor decomposition, which helps uncover latent structures or estimate missing values in data. Several tensor decomposition techniques have been developed, each characterized by a distinct way of decomposing tensor components. Common methods include canonical polyadic (CP) decomposition \cite{RM24,hitchcock1927}, tucker decomposition \cite{tucker1966mathematical}, tensor train decomposition \cite{oseledets2011tensor}, tensor ring decomposition \cite{yuan2019tensor}, and block term decomposition (BTD) \cite{rontogiannis2021block}.  Among these, CP decomposition is one of the most fundamental and widely used approaches, see, e.g. \cite{ CM,goulart2015tensor,TG,RM22,XL,ND}.
 It represents a tensor as a sum of rank-1 components, making it a natural extension of matrix factorization to higher-order data. Due to its strong interpretability and broad applicability, CP decomposition has been extensively studied, leading to various numerical algorithms for its computation. The most commonly used algorithm is the alternating least squares (ALS) \cite{RM23,RM24} method, which iteratively updates factor matrices to minimize reconstruction error \cite{Tb,Tl}.
While the ALS algorithm is simple and easy to implement, its convergence speed is relatively slow, which limits its practical efficiency. To address this issue, various improvements to the CP-ALS algorithm have been proposed to enhance both efficiency and convergence. For example, Chen et al. \cite{Chen2011} introduced a higher-order search direction and an adaptive step size selection strategy to accelerate the ALS algorithm. Li et al. \cite{Li2013} developed a nonlinear least squares method based on Gauss-Newton optimization to minimize the $L_{2}$ norm of the approximation residual, yielding a more effective low-rank CP decomposition. Rajih et al. \cite{Rajih2008} extended the ALS algorithm by incorporating a more sophisticated extrapolation technique that leverages nonlinear parameter trends and updates all parameter sets simultaneously. 
 In addition to the ALS algorithm, other algorithms have also been developed for solving CP decomposition. Notable examples include those in \cite{comon2009tensor, Navasca2008, Tomasi2006,tomasi2006comparison}.
 
Specifically, algorithms derived from the ALS framework are required to solve a series of quadratic subproblems, which often involve matrix inversion to obtain the solutions. However, explicit matrix inversion is generally discouraged due to two main drawbacks. First, it incurs a significant computational cost, often scaling at $O(n^3)$ for dense matrices, which becomes prohibitively expensive as the matrix size increases. Second, it is more susceptible to numerical instability, especially when applied to ill-conditioned matrices, thereby compromising the reliability of the solution.  To mitigate these issues, researchers have explored alternative strategies to improve the efficiency and stability of ALS-based methods. One effective approach is to replace the conventional normal equation method with QR decomposition when solving quadratic subproblems \cite{RM20}. Unlike the normal equation method, which requires computing the matrix inverse explicitly, QR decomposition leverages the orthogonality property $\mathbf{Q}^T\mathbf{Q}=\mathbf{I}$ of orthogonal matrices. This not only avoids the numerical instability caused by matrix inversion, but also improves the overall robustness of the algorithm. By integrating QR decomposition into the ALS framework, the solution process becomes more stable and computationally efficient, paving the way for more reliable tensor decomposition methods.

\subsection{Motivation and contribution}
The CP decomposition employing QR decomposition principles, denoted as CP-ALS-QR \cite{RM20}, solves for the factor matrices by iteratively solving the linear equations of the upper triangular matrix, specifically using the relation $\mathbf{R}\mathbf{X}=\mathbf{Q}^T\mathbf{B}$, where $\mathbf{X}$ represents the factor matrix and $\mathbf{B}$ corresponds to the unfolded tensor. This formulation leverages the triangular structure of $\mathbf{R}$, transforming the solution process into a sequence of linear equations, which not only greatly reduces the computational burden of matrix inversion but also endows the CP-ALS-QR framework with better numerical stability compared to the CP-ALS framework. However, the CP-ALS-QR algorithm requires multiple tensor-times-matrix (TTM) operations during the iteration process, significantly increasing the computational overhead. This challenge becomes particularly prominent in high-rank scenarios, where its computational complexity markedly exceeds that of the traditional CP-ALS algorithm.

To overcome this limitation, we investigated the use of dimension tree optimization techniques for tensor TTM operations. Although the dimension tree has shown improvements in computational efficiency and complexity both theoretically and experimentally, especially in synthetic tensor scenarios, their performance can fluctuate when applied to real-world tensor datasets. Furthermore, during iterative updates using the dimension tree, some intermediate tensors are underutilized, resulting in inefficient use of the dimension tree. On the other hand, inspired by the widespread use of acceleration techniques, such as Nesterov acceleration \cite{Nesterov1983}, momentum acceleration \cite{evans2023blockwise}, and extrapolation acceleration \cite{Ang2019}, we apply some acceleration techniques to speed up the CP-ALS-QR algorithm. However, considering the inherent computational complexity of the algorithm, we mainly explore simpler forms of acceleration to ensure that the computational complexity of the algorithm does not increase significantly, or can even be neglected. 

The key contributions of this paper can be summarized as follows.

\begin{itemize}
\item We propose the restructured dimension tree to enhance the reuse of intermediate tensors within the dimension tree. When integrated into the CP-ALS-QR algorithm, it significantly improves iteration efficiency by increasing the utilization of intermediate tensors. Theoretical analysis in Subsection \ref{Cost} further demonstrates that, compared with the standard dimension tree, the restructured dimension tree reduces computational complexity. Compared with the standard dimension tree, the restructured dimension tree can reduce the major portion of the computational cost incurred by the TTM operations of tensors of any order by 33\%.

\item We propose an efficient and concise customized extrapolation for the matrix $\mathbf{Q}_0$ within the CP-ALS-QR algorithm. This extrapolation technique leverages the unique structural properties of $\mathbf{Q}_0$, where the first element in its first row and first column is 1, and all other elements in the first row and first column are 0. It is specifically used to accelerate the extrapolation method of the CP-ALS-QR algorithm.

\item  We introduce the algorithm ALS-QR-BRE (see Algorithm \ref{ALS-QR-BRE}), which combines the restructured dimension tree with customized extrapolation techniques. Experimental evaluations on both synthetic and real-world tensor datasets indicate that ALS-QR-BRE outperforms CP-ALS-QR in terms of both iteration speed and fitting accuracy.  In experiments on real-world datasets, its iteration speed is improved by approximately one to two times.

\end{itemize}

\section{Background}\label{sect-2}
In this section, we introduce the essential mathematical notations and fundamental concepts used throughout the paper. We first define the notation for tensors and their basic operations, followed by a discussion of key properties of the  CP decomposition. Additionally, we summarize two relevant CP decomposition algorithms: CP-ALS \cite{TG} and CP-ALS-QR \cite{RM20}.
\vspace{-2pt}
\subsection{Notation and definitions}
The analysis in this paper employs matrix multiplication, tensor algebra, and specialized notation for tensor operations, adhering to the symbolic conventions established in \cite{TG}. We denote scalars by lowercase letters (e.g., a), vectors by boldface lowercase letters (e.g., $\mathbf{a}$), matrices by boldface uppercase letters (e.g., $\mathbf{A}$), and tensors by calligraphic uppercase letters (e.g., $\mathcal{X}$).

The kronecker product of two matrices $\mathbf{A} \in \mathbb{R}^{I \times J}$ and $\mathbf{B} \in \mathbb{R}^{K \times L}$ is denoted by $\mathbf{A} \otimes \mathbf{B} \in \mathbb{R}^{(IK) \times (JL)}$, with the entries given by $\left[ \mathbf{A} \otimes \mathbf{B} \right]_{(K(i-1)+k, L(j-1)+l)} = \mathbf{A}(i,j)\mathbf{B} (k,l)$. For matrices $\mathbf{A} \in \mathbb{R}^{I \times J}$ and $\mathbf{B} \in \mathbb{R}^{K \times J}$, their khatri-rao product results in a matrix of size $(IK) \times J$, defined as $\mathbf{A} \odot \mathbf{B} = [\mathbf{a}_1 \otimes \mathbf{b}_1, \dots, \mathbf{a}_J \otimes \mathbf{b}_J]$. For any matrices $\mathbf{A} \in \mathbb{R}^{J \times I}$, $\mathbf{B} \in \mathbb{R}^{K \times I}$, $\mathbf{C} \in \mathbb{R}^{K \times J}$, and $\mathbf{D} \in \mathbb{R}^{J \times K}$, within the context of kronecker and khatri-rao product operations,
the following identity holds:
\vspace{-5pt}
\begin{equation}\label{eq-1}
    (\mathbf{C} \otimes \mathbf{D})(\mathbf{A} \odot \mathbf{B} ) = (\mathbf{C}\mathbf{A})\odot (\mathbf{D}\mathbf{B}).
\end{equation}
The hadamard product of two matrices $\mathbf{A} \in \mathbb{R}^{I \times J}$ and $\mathbf{B} \in \mathbb{R}^{I \times J}$ is the elementwise product, denoted as $\mathbf{A} \ast \mathbf{B} \in \mathbb{R}^{I \times J}$, where $[\mathbf{A} \ast \mathbf{B}]{(i,j)} = \mathbf{A}{(i,j)} \mathbf{B}{(i,j)}$. The mode-$n$ product of a tensor $\mathcal{X} \in \mathbb{R}^{I_{1} \times \cdots \times I_{N}}$ and a matrix $\mathbf{B} \in \mathbb{R}^{J \times I_{n}}$ is defined by $\mathcal{Y} = \mathcal{X} \times_{n} \mathbf{B} \in \mathbb{R}^{I_{1} \times \cdots \times I_{n-1} \times J \times I_{n+1} \times \cdots \times I_{N}}$, with entries given by $\mathcal{Y}(i_{1}, \dots, j, \dots, i_{N}) = \sum_{i_{n}=1}^{I_{n}} \mathcal{X}(i_{1}, \dots, i_{n}, \dots, i_{N}) \mathbf{B}(j, i_{n})$. This operation, commonly referred to as  tensor-times-matrix (TTM) operation, extends naturally to scenarios where an $N$-mode tensor is multiplied by multiple matrices along different modes. This extended operation, known as Multi-TTM, is formally defined as: 
\begin{equation}
    \mathcal{Y} = \mathcal{X} \times_{1}\mathbf {B}_1 \times_{2} \mathbf{B}_2 \times_{3} \cdots \times_{N}\mathbf {B}_N,
\end{equation}
where $\mathcal{X} \in \mathbb{R}^{I_{1}\times \cdots \times I_{N}}$, $\mathcal{Y} \in \mathbb{R}^{R_{1}\times \cdots \times R_{N}}$ and $\mathbf{B_i}\in\mathbb{R}^{R_{i}\times I_{i}} $. This is expressed in the matricized form as:
\begin{equation}\label{eq-3}
    \mathbf{Y}_{(n)}= \mathbf{X}_{(n)}(\mathbf{B}_{N}\otimes \cdots\otimes \mathbf{B}_{n-1}\otimes \mathbf{B}_{n+1}\otimes \cdots\otimes \mathbf{B}_{1})^{T}, 
\end{equation}
where $\mathbf{X}_{(n)}$ and $\mathbf{Y}_{(n)}$ represent the mode-$n$ unfolding of tensors $\mathcal{X}$ and $\mathcal{Y}$, respectively.

\subsection{CP decomposition}
The CP decomposition \cite{RM24,hitchcock1927} is a higher-order extension of matrix singular value decomposition. It expresses a higher-order tensor as a sum of a finite number of rank-1 tensors. For an $N$-order tensor $\mathcal{X} \in \mathbb{R}^{I_{1} \times \cdots \times I_{N}}$, the rank-$R$ CP decomposition of $\mathcal{X}$ is an approximation
\vspace{-3pt}
\begin{equation} 
  \mathcal{X}\approx \sum_{r=1}^{R} \mathbf{\lambda}_{r} \mathbf{a}_{r}^{(1)}\circ  \mathbf{a}_{r}^{(2)} \circ \cdots \circ \mathbf{a}_{r}^{(N)},
  \vspace{-3pt}
\end{equation}
where $\mathbf{a}_{r}^{(n)}\in \mathbb{R}^{I_{n}}$ are unit vectors with weight $\lambda_r$, and $\circ$ denotes the outer product. The matrix formed by all the $\mathbf{a}_{r}^{(n)}$ vectors is called the factor matrix, that is , $A^{(n)} = \left [\mathbf{a}_{1}^{(n)} \quad  \mathbf{a}_{2}^{(n)} \quad \cdots  \quad \mathbf{a}_{r}^{(n)}   \right ] $.
So the CP decomposition  can also be denoted by
\begin{equation}
    \mathcal{X}\approx \left [ \left [ \lambda :\mathbf{A}^{(1)}, \mathbf{A}^{(2)},\cdots, \mathbf{A}^{(N)}  \right ]  \right ] ,
    \vspace{-5pt}
\end{equation}
where $\lambda  \in \mathbb{R}^{R}$ denotes the weight vector.

\subsection{CP-ALS}
The ALS \cite{RM23, RM24} is the most widely used CP decomposition algorithm. The core idea of this algorithm is to update one factor matrix at a time while keeping the others fixed. For the 
$n$th factor matrix, the ALS algorithm iteratively updates by solving the following quadratic subproblem:
\begin{equation}\label{eq-6}
\operatorname*{min}_{\hat{\mathbf{A}}^{(n)}} \left \|  \mathbf{X}_{(n)}- \hat{\mathbf{A}}^{(n)}\mathbf{P}^{(n)T}\right \|_F,
\end{equation}
where $\hat{\mathbf{A}}^{(n)} = \mathbf{A}^{(n)} \cdot \text{diag}(\lambda)$ and $\mathbf{X}_{(n)}$ denotes the mode-$n$ unfolding of the tensor $\mathcal{X} \in \mathbb{R}^{I_{1} \times \cdots \times I_{N}}$. The  $\mathbf{P}^{(n)} \in \mathbb{R}^{W_{n} \times R}$, where $W_{n} = I_{1} \times \cdots \times I_{n-1} \times I_{n+1} \times \cdots \times I_{N}$ and $R$ is the approximative rank of the $\mathcal{X}$, is formed by the Khatri-Rao product of the other factor matrices.
\begin{equation}
\mathbf{P}^{(n)}=\mathbf{A}^{(1)} \odot \cdots \odot \mathbf{A}^{(n-1)} \odot \mathbf{A}^{(n+1)} \odot  \cdots \odot \mathbf{A}^{(N)}.
\end{equation}
The subproblem (\ref{eq-6}) can be transformed into the following form:
\vspace{-3.25pt}
\begin{equation}\label{eq_8}
    \hat{\mathbf{A}}^{(n)}\mathbf{P}^{(n)T} \cong  \mathbf{X}_{(n)}.
\vspace{-3.25pt}
\end{equation}
To compute $\hat{\mathbf{A}}^{(n)}$, we need to solve the linear system defined in (\ref{eq_8}). However, since $\mathbf{P}^{(n)T}$ is not a square matrix, it does not have a direct inverse. To address this, we pre-multiply both sides of (\ref{eq_8}) by $\mathbf{P}^{(n)}$, forming the normal equation. This results in a symmetric positive semi-definite system involving the square matrix $\Gamma^{(n)}$, which can be efficiently inverted if it is well-conditioned. Consequently, by (\ref{eq_8}) , we derive
\vspace{3pt}
\begin{equation}\label{eq-9}
 \hat{\mathbf{A}}^{(n)}\Gamma^{(n)} =  \mathbf{X}_{(n)}\mathbf{P}^{(n)},
\end{equation}
where $\Gamma^{(n)} \in \mathbb{R}^{R \times R}$ can be computed via
\begin{equation}
\Gamma^{(n)}=\mathbf{S}_{1}\ast \cdots \ast \mathbf{S}_{n-1}\ast \mathbf{S}_{n+1}\ast \cdots \ast \mathbf{S}_{N},
\end{equation}
where $\mathbf{S}_i=(\mathbf{A}^{(i)})^T\mathbf{A}^{(i)},i=1,\cdots,n-1,n+1,\cdots,N$.

By multiplying both sides of (\ref{eq-9}) by $\Gamma^{(n)-1}$, we obtain $\hat{\mathbf{A}}^{(n)}$. Subsequently, normalizing the columns of $\hat{\mathbf{A}}^{(n)}$ yields the factor matrix $\mathbf{A}^{(n)}$. The detailed implementation of the ALS algorithm is provided in Algorithm \ref{CP-ALS}.

\begin{algorithm}
    \caption{CP-ALS}
{\bfseries Input:} Tensor $\mathcal{X}\in \mathbb{R}^{I_{1}\times \cdots \times I_{N}}$, rank $R$, stopping criteria $\epsilon$, iteration count $m$.\\
{\bfseries Initialize:} Factor matrices $\mathbf{A}^{(1)}, \cdots, \mathbf{A}^{(N)}$, gram matrices $\mathbf{S}_{1}=(\mathbf{A}^{(1)})^{T} \mathbf{A}^{(1)},\cdots,\mathbf{S}_{N}=(\mathbf{A}^{(N)})^{T}\mathbf{A}^{(N)}$.
\begin{algorithmic}[1]
\For {$k=0,1,\dots m$} 
\For {$n=1,2,\dots N$}
\State  $\Gamma^{(n)}=\mathbf{S}_{1}* \cdots* \mathbf{S}_{n-1}* \mathbf{S}_{n+1}* \cdots * \mathbf{S}_{N}$.
\State $\mathbf{P}^{(n)}=\mathbf{A}^{(1)}\odot \cdots \odot \mathbf{A}^{(n-1)}\odot \mathbf{A}^{(n+1)}\odot\cdots\odot \mathbf{A}^{(N)}$.
\State $\mathbf{M}_{n}=\mathbf{X}_{(n)}\mathbf{P}^{(n)}$.
\State Solve $\hat{\mathbf{A}}^{(n)}\Gamma^{(n)}=\mathbf{M}_{n}$ to obtain $\hat{\mathbf{A}}^{(n)}$.
\State Normalize columns of $\hat{\mathbf{A}}^{(n)}$ to obtain factor matrix $\mathbf{A}^{(n)}$.
\State Recompute gram matrix $\mathbf{S}_{n}=(\mathbf{A}^{(n)})^{T}\mathbf{A}^{(n)}$ for updated factor matrix $\mathbf{A}^{(n)}$.
\EndFor
\State Compute the $fitness$.
\If{$fitness\ge \epsilon$}
\State Break
\EndIf
\EndFor
\end{algorithmic}
{\bfseries Output:} Matrices $\{ \mathbf{A}^{(1)} \cdots \mathbf{A}^{(N)}\}$.
\label{CP-ALS}
\end{algorithm}

\subsection{CP-ALS-QR} 
The CP-ALS algorithm solves the linear system in (\ref{eq_8}) using the normal equation, which involves forming $\Gamma^{(n)} = \mathbf{P}^{(n)T} \mathbf{P}^{(n)} $ and explicitly computing its inverse. In contrast,  {the CP-ALS-QR  \cite{RM20}} improves numerical stability by replacing the normal equation approach with QR decomposition, thereby avoiding direct matrix inversion.

The transformation of $\mathbf{P}^{(n)T}$ is carried out as follows: First, we compute the compact QR factorization of each individual factor matrix, such that $\mathbf{A}^{(i)} = \mathbf{Q}_{i} \mathbf{R}_{i}$. Consequently, $\mathbf{P}^{(n)}$ can be reformulated in the following manner:
\begin{equation}
\begin{aligned}
   \mathbf{P}^{(n)}&=\mathbf{A}^{(N)}\odot \cdots \odot \mathbf{A}^{(n+1)}\odot \mathbf{A}^{(n-1)}\odot \cdots \odot \mathbf{A}^{(1)}
        \\&=\mathbf{Q}_{N}\mathbf{R}_{N}\odot \cdots \odot \mathbf{Q}_{n+1}\mathbf{R}_{n+1}\odot \mathbf{Q}_{n-1}\mathbf{R}_{n-1}\odot \cdots \odot \mathbf{Q}_{1}\mathbf{R}_{1}
        \\&=(\mathbf{Q}_{N}\otimes \cdots \otimes \mathbf{Q}_{n+1}\otimes \mathbf{Q}_{n-1}\otimes \cdots \otimes \mathbf{Q}_{1})\underbrace{ (\mathbf{R}_{N}\odot \cdots \odot \mathbf{R}_{n+1}\odot \mathbf{R}_{n-1}\odot \cdots \odot \mathbf{R}_{1})}_{\mathbf{Z}_{n}},
\end{aligned}
\end{equation}
where the last equality follows from (\ref{eq-1}). Second, we compute the QR decomposition of the matrix $\mathbf{Z}_{n}=\mathbf{Q}_{0}\mathbf{R}_{0}$. Thus, we have
\vspace{-2pt}
\begin{equation}
\begin{aligned}\label{eq-11}
\mathbf{P}^{(n)}=&(\mathbf{Q}_{N}\otimes \cdots \otimes \mathbf{Q}_{n-1}\otimes \mathbf{Q}_{n+1}\otimes \cdots \otimes \mathbf{Q}_{1})
        (\mathbf{R}_{1}\odot \cdots \odot \mathbf{R}_{n-1}\odot \mathbf{R}_{n+1}\odot \cdots \odot \mathbf{R}_{N})
      \\=&\underbrace{(\mathbf{Q}_{N}\otimes \cdots \otimes \mathbf{Q}_{n+1}\otimes \mathbf{Q}_{n-1}\otimes \cdots \otimes \mathbf{Q}_{1})\mathbf{Q}_{0}}_{\mathbf{Q}} \underbrace{\mathbf{R}_{0}}_{\mathbf{R}},
\end{aligned} 
\end{equation}
where $\mathbf{Q}$ has orthonormal columns, and $\mathbf{R}$ is upper triangular.

By applying the aforementioned QR decomposition to $\mathbf{P}^{(n)T}$, linear system (\ref{eq_8}) can be reformulated as follows:
\begin{equation}\label{eq-13}
    \hat{\mathbf{A}}^{(n)}\mathbf{R}^T\mathbf{Q}^T =  \mathbf{X}_{(n)},
\end{equation}
where $\mathbf{Q}={(\mathbf{Q}_{N}\otimes \cdots \otimes \mathbf{Q}_{n+1}\otimes \mathbf{Q}_{n-1}\otimes \cdots \otimes \mathbf{Q}_{1})\mathbf{Q}_{0}}$, $\mathbf{R}=\mathbf{R}_0$. Given that the matrix $\mathbf{Q}$ is orthonormal, it follows that $\mathbf{Q}^T\mathbf{Q}=\mathbf{I}$. Consequently, by multiplying both sides of (\ref{eq-13}) by $\mathbf{Q}$, we obtain 
\begin{equation}\label{eq-14}
    \hat{\mathbf{A}}^{(n)}\mathbf{R}^T = \mathbf{Y}_{(n)}\mathbf{Q},
\end{equation}
where $\mathbf{Y}_{(n)}= \mathbf{X}_{(n)}(\mathbf{Q}_{N}\otimes \cdots \otimes \mathbf{Q}_{n-1}\otimes \mathbf{Q}_{n+1}\otimes \cdots \otimes \mathbf{Q}_{1})$. The $\mathbf{Y}_{(n)}$ is the mode-$n$ unfolding of $\mathcal{Y} =\mathcal{X}\times_{1}\mathbf{Q}_{1}^{T}\times_{2}\cdots \times_{i-1}\mathbf{Q}_{i-1}^{T}\times_{i+1}\mathbf{Q}_{i+1}^{T}\times_{n+2}\cdots \times_{N}\mathbf{Q}_{N}^{T}$. The detailed implementation of the ALS algorithm is provided in Algorithm \ref{CP-ALS-QR}.

\begin{algorithm}
\caption{CP-ALS-QR}
{\bfseries Input:} Tensor $\mathcal{X}\in \mathbb{R}^{I_{1}\times \cdots \times I_{N}}$, rank $R$, stopping criteria $\epsilon$, iteration count $m$.\\
{\bfseries Initialize:} Factor matrices $\mathbf{A}^{(1)} \cdots \mathbf{A}^{(N)}$, calculate QR-decomposition $\mathbf{Q}_{1}\mathbf{R}_{1}\cdots \mathbf{Q}_{N}\mathbf{R}_{N}$ of factor matrices.
\begin{algorithmic}[1]
\For {$k=0,1,\dots m$}   
\For {$n=1,2,\dots N$}
 \State  $\mathbf{Z}_{n}=\mathbf{R}_{N}\odot \cdots\odot \mathbf{R}_{n+1}\odot \mathbf{R}_{n-1}\odot \cdots \odot \mathbf{R}_{1}$.
 \State  Calculate QR decomposition $\mathbf{Z}_{n}=\mathbf{Q}_{0}\mathbf{R}_{0}$.
\State $\mathcal{Y} =\mathcal{X}\times_{1}\mathbf{Q}_{1}^{T}\times_{2}\cdots \times_{n-1}\mathbf{Q}_{n-1}^{T}\times_{n+1}\mathbf{Q}_{n+1}^{T}\times_{n+2}\cdots \times_{N}\mathbf{Q}_{N}^{T}$.
\State  $\mathbf{V}_{n}=\mathbf{Y}_{(n)}\mathbf{Q}_{0}$.
\State  Solve $\hat{\mathbf{A}}^{(n)}\mathbf{R}_{0}^{T}=\mathbf{V}_{n}$ to obtain $\hat{\mathbf{A}}^{(n)}$.
\State  Normalize columns of $\hat{\mathbf{A}}^{(n)}$ to obtain factor matrix $\mathbf{A}^{(n)}$.
\State Recompute QR-decomposition for updated factor matrix $\mathbf{A}^{(n)}=\mathbf{Q}_{n}\mathbf{R}_{n}$.    
\EndFor
\State Compute the $fitness$.
\If{$fitness\ge \epsilon$}
\State Break
\EndIf
\EndFor
\end{algorithmic}
{\bfseries Output:} Matrices $\{{\mathbf{A}^{(1)} \cdots \mathbf{A}^{(N)}}\}$.
\label{CP-ALS-QR}
\end{algorithm}

\section{Restructured dimension tree and  customized extrapolation }\label{sect-3}
In this section, we introduce the restructured dimension tree and the customized extrapolation design for matrix $\mathbf{Q}_0$ in line 4 of Algorithm \ref{CP-ALS-QR}. Based on these techniques, we develop an ALS-QR-BRE algorithm for the CP decomposition. Additionally, we give a complexity analysis of the restructured dimension tree.

\subsection{Restructured dimension tree} \label{BR_subsection}
\begin{figure*}[!htbp]
\centering  
\subfloat{
\label{Fig.sub-1}
\includegraphics[width=4.5cm,height = 4.5cm]{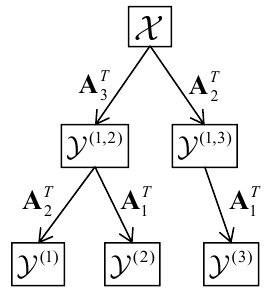}}
\subfloat{
\label{Fig.sub.2}
\includegraphics[width=6cm,height = 5.5cm]{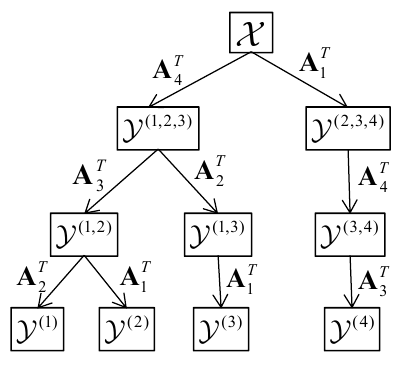}}
\vspace{-3pt}
\caption{The canonical dimension tree architecture of third-order and fourth-order tensors.}
\label{fig-1}
\end{figure*}

The dimension tree provides a hierarchical partitioning of mode indices in an $N$-dimensional tensor, enabling more efficient computation in tensor decomposition. Originally, it was introduced to facilitate hierarchical tucker decomposition \cite{Grasedyck2010}, which offers a structured way to represent tensors while reducing computational complexity. Compared to the conventional tucker decomposition, the dimension tree approach is particularly advantageous for high-order tensors, where direct tucker decomposition becomes impractical due to its high computational cost. A key component of the dimension tree structure is the intermediate tensors, which are also called TTMc tensors. These TTMc tensors, denoted as $\mathcal{Y}^{(i_{1},i_{2},...,i_{m})}$, are defined as follows:

\begin{equation}\label{eq-12}
\mathcal{Y}^{(i_{1},i_{2},...,i_{m})}=\mathcal{X}\times_{n\in{\left \{ 1,2,...,N \right \} \setminus \left \{i_{1},i_{2},...,i_{m}  \right \} }}\mathbf{A}_{n}^{T},
\end{equation}
where $\mathcal{X}$ is contracted with all the matrices $\mathbf{A}_{n}$ except $\mathbf{A}_{i_{1}},...,\mathbf{A}_{i_{m}}$. The $\mathbf{A}_{n}$ is with respect to $\mathbf{Q}_{n}$ in line 5 of Algorithm \ref{CP-ALS-QR}.
Building upon the definition of intermediate tensors delineated in (\ref{eq-12}), the  dimension trees for third-order and fourth-order tensors are depicted schematically in Figure \ref{fig-1}.

\begin{figure}[!h]
\centering  
\subfloat{
\includegraphics[scale=0.80]{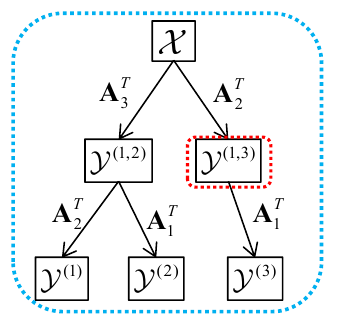}}
\subfloat{
\includegraphics[scale=0.80]{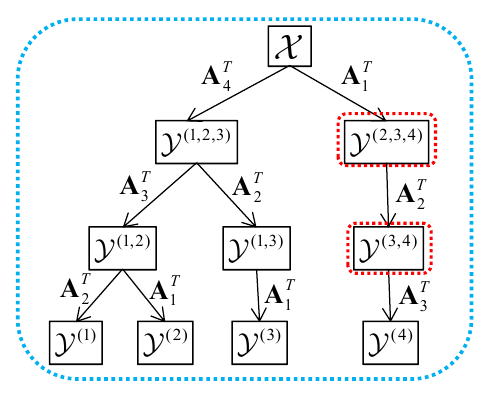}}
\caption{The tensor that remains underutilized is emphasized within the red-dashed box. (a) Shows the underutilized tensor $\mathcal{Y}^{(1,3)}$. (b) Shows the underutilized tensors $\mathcal{Y}^{(3,4)}$ and $\mathcal{Y}^{(2,3,4)}$.}
\label{fig-2}
\end{figure}

\vspace{-10pt}
As shown in Figure \ref{fig-1}, the dimension tree structure does not fully exploit certain intermediate tensors during the iteration process, leading to inefficiencies. For instance, in Figure \ref{fig-1}(a), the tensor $\mathcal{Y}^{(1,3)}$ remains unexploited to updating the factor matrix $\mathbf{A}^{(1)}$, which results in unnecessary memory usage. To address this issue, we propose the restructured dimension tree to better utilize intermediate tensors and enhance TTM computational efficiency. The key principles of this approach are maximizing the reuse of intermediate tensors, ensuring that previously computed tensors are effectively utilized in subsequent iterations, and reducing reliance on the original tensor by using intermediate tensors to update factor matrices associated with different indices, minimizing the need to repeatedly access the original tensor. Figure \ref{fig-2} illustrates which intermediate tensors are not fully utilized during the iteration process.

\begin{figure}[!htbp]
   \centering 
   \includegraphics[scale=0.8]{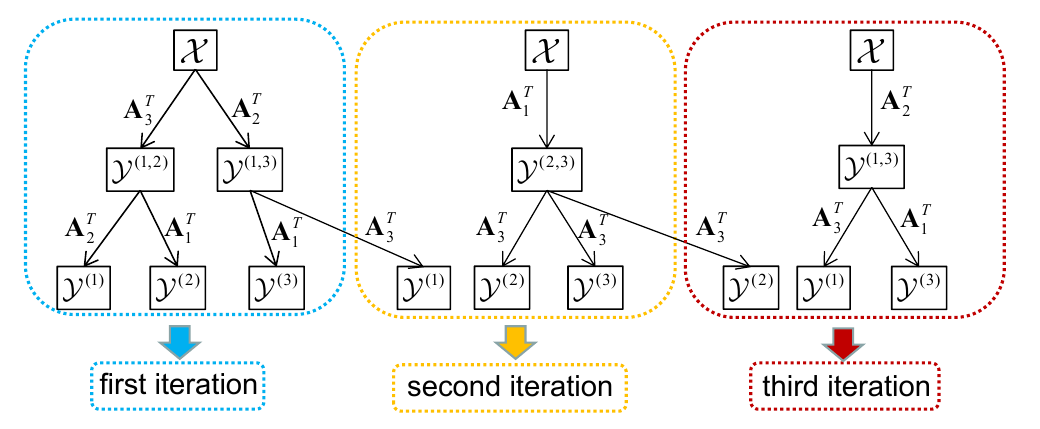}
   \caption{The iterative structure of  restructured dimension tree for third-order tensor. }
   \label{fig-3}
\end{figure}

In the case of third-order tensors, unlike the traditional dimension tree, which repeatedly applies the same structural approach, the proposed restructured dimension tree starts the second iteration from the underutilized tensor $\mathcal{Y}^{(1,3)}$. By the third iteration, all tensors have been fully utilized, eliminating the need to access underutilized tensors. Consequently, rather than restarting from the original tensor, the iteration proceeds directly from the intermediate tensor $\mathcal{Y}^{(2,3)}$. If the third iteration were to begin from the original tensor, an additional TTM computation would be required every two iterations, compared to starting from an intermediate tensor. Therefore, to minimize computational overhead, it is more efficient to avoid restarting iterations from the original tensor. Figure \ref{fig-3} illustrates the restructured dimension tree for third-order tensor.

\begin{figure}[!htbp]
   \centering 
   \includegraphics[scale=0.70]{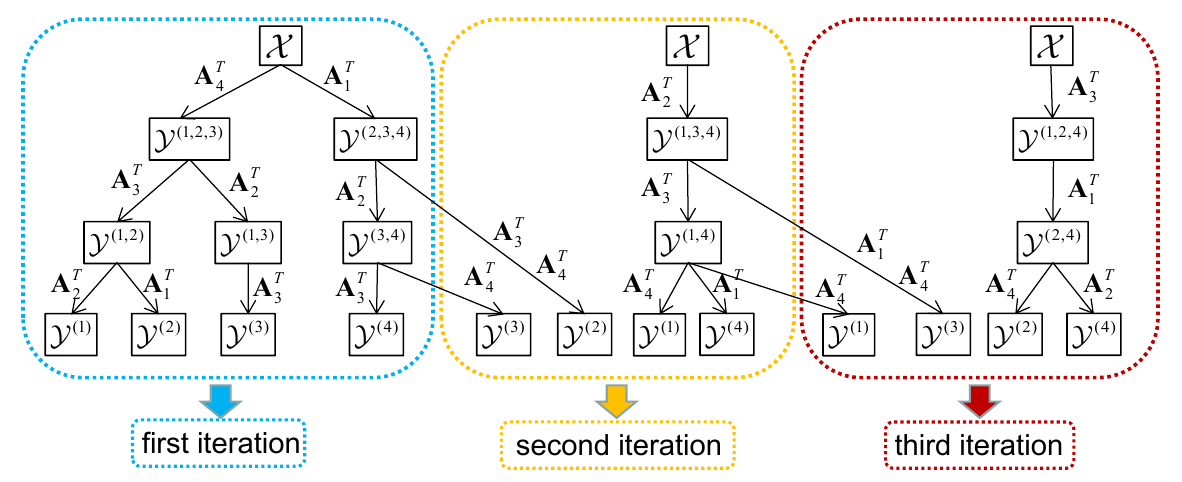}
   \caption{The iterative structure of restructured dimension tree for fourth-order tensor. }
   \label{fig-4}
\end{figure}

For fourth-order tensors, the restructured dimension tree updates iteratively by combining intermediate tensors with underutilized ones. In each iteration, we prioritize reusing underutilized tensors to compute the corresponding factor matrices. In the second iteration, the intermediate tensor $\mathcal{Y}^{(3,4)}$, generated in the first iteration, remains unused for updating $\mathbf{A}^{(3)}$, so the second iteration begins from $\mathcal{Y}^{(3,4)}$. By the third iteration, the intermediate tensor $\mathcal{Y}^{(1,3,4)}$ from the second iteration is still underutilized and not yet applied to compute $\mathbf{A}^{(3)}$. Therefore, the third iteration begins by utilizing $\mathcal{Y}^{(1,4)}$ to update $\mathbf{A}^{(1)}$, followed by the use of $\mathcal{Y}^{(1,3,4)}$ to compute $\mathbf{A}^{(3)}$. Figure \ref{fig-4} illustrates the refined iterative structure for the restructured dimension tree in the fourth-order tensor.

In this subsection, we present a comprehensive depiction of the restructured dimension tree. The core principle driving the construction of this method is based on alternating iterative updates of the factor matrices, wherein the restructured dimension tree systematically performs alternating updates of the factor matrices on the intermediate tensors.

\subsection{Customized extrapolation} \label{extra_subsection}

Extrapolation acceleration techniques \cite{Nesterov1983} are well-studied methods for improving the efficiency of iterative optimization algorithms and have been extensively used in the field of convex optimization. In recent years, significant progress has been made in incorporating extrapolation strategies into the CP-ALS algorithm. Notable approaches include line search methods \cite{Bro1997}, momentum-based techniques \cite{Mitchell2020}, and block momentum acceleration \cite{evans2023blockwise}. These methods accelerate convergence by leveraging information from previous iterations to refine the current solution, thereby enhancing the overall optimization process.

We also apply extrapolation acceleration to improve the convergence of the algorithm. Nevertheless, considering the computational overhead of the original method, we confine our approach to a basic extrapolation strategy. Therefore, we propose introducing an extrapolation strategy into Algorithm \ref{CP-ALS-QR} to improve numerical efficiency, which is implemented in line 4 of Algorithm \ref{CP-ALS-QR}, as follows
\begin{equation}\label{eq-18} 
\hat{\mathbf{Q}}_0^k = \mathbf{Q}_0^k + \beta_k ( \mathbf{Q}_0^k - \mathbf{Q}_0^{k-1} ). 
\end{equation} 
However, empirical analysis shows that the effectiveness of this extrapolation mechanism in improving tensor approximation accuracy is restricted by the inherent structure of $\mathbf{Q}_0$. In particular, $\mathbf{Q}_0$ exhibits a special form in which the entry at the first row and first column equals 1, while all other entries in the first row and first column are 0. For instance, in the case of a $4 \times 4$ matrix, its form is as follows
\vspace{5pt}
$$
Q_0=\begin{bmatrix}
  1& 0 &0  &0 \\
  0& q_{22} & q_{23} &q_{24} \\
  0&  0& q_{33} &q_{34} \\
  0&  0&  0&q_{44}
\end{bmatrix}.
\vspace{5pt}
$$

\noindent Therefore, the extrapolation operation (\ref{eq-18}) preserves the first row and first column of $\mathbf{Q}_0$, which significantly limits the effectiveness of extrapolation. To overcome this limitation, we propose a customized extrapolation approach tailored for $\mathbf{Q}_0$, which better accommodates its structure while keeping the additional computational cost of the original algorithm nearly negligible. The customized extrapolation for $\mathbf{Q}_0$ is as follows:
\vspace{5pt}
\begin{equation}\label{e_19}
\hat{\mathbf{Q}}_0^k = \mathbf{Q}_0^k + \beta ( \mathbf{Q}_0^k - \alpha \mathbf{Q}_0^{k-1} ), 
\end{equation} 
where $\beta$ and $\alpha$ represent tunable hyperparameters. Unlike the traditional momentum extrapolation (\ref{eq-18}), the modified form (\ref{e_19}) slightly reduces acceleration but enhances numerical stability and better aligns with the structural constraints of $\mathbf{Q}_0$.

\begin{remark}
 The proposed extrapolation is applied to a single variable and is specifically tailored to the underlying data structure of that variable.
\end{remark}

\subsection{ALS-QR-BRE}
Through integrating the optimized dimension tree architecture (see Subsection \ref{BR_subsection}) and the customized extrapolation strategy (see Subsection \ref{extra_subsection}), we propose a novel algorithm, ALS-QR-BRE. This algorithm is based on two key principles: (i) the restructured dimension tree accelerates computations in line 5 of Algorithm \ref{CP-ALS-QR}, leading to improved computational efficiency; (ii) the customized extrapolation technique is applied to matrix $\mathbf{Q}_0$ in line 4 of Algorithm \ref{CP-ALS-QR} to enhance approximation accuracy. The full implementation is provided in Algorithm \ref{ALS-QR-BRE}.

\begin{algorithm}[!ht]\label{ALS-QR-BR}
\caption{ALS-QR-BRE}
{\bfseries Input:} Tensor $\mathcal{X}\in \mathbb{R}^{I_{1}\times \cdots \times I_{N}}$, rank $R$, stopping criteria $\epsilon$, iteration count $m$.\\
{\bfseries Initialize:} Factor matrices $\mathbf{A}^{(1)} \cdots \mathbf{A}^{(N)}$, calculate QR-decomposition $\mathbf{Q}_{1}\mathbf{R}_{1},\cdots, \mathbf{Q}_{N}\mathbf{R}_{N}$ of factor matrices.
\begin{algorithmic}[1]
\For {$k=0,1,\dots, m$}   
\For {$n=1,2,\dots, N$}
 \State  $\mathbf{Z}_{n}=\mathbf{R}_{N}\odot \cdots\odot \mathbf{R}_{n+1}\odot \mathbf{R}_{n-1}\odot \cdots \odot \mathbf{R}_{1}$.
 \State  Calculate QR-decomposition $\mathbf{Z}_{n}=\mathbf{Q}_{0}^k\mathbf{R}_{0}$.
 \State  $\hat{\mathbf{Q}}_0^k=\mathbf{Q}_0^k + \beta ( \mathbf{Q}_0^k - \alpha \mathbf{Q}_0^{k-1} )$.  {\color{gray}$\hspace{4.5cm}\,\leftarrow$ Subsection \ref{extra_subsection}}
\State   Calculate $\mathcal{Y} $ using the \textbf{Restructured dimension tree}.  {\color{gray} $\quad \leftarrow$ Subsection \ref{BR_subsection}}
\State  $\mathbf{V}_{n}=\mathbf{Y}_{(n)}\hat{\mathbf{Q}}_{0}^k$.
\State  Solve $\hat{\mathbf{A}}^{(n)}\mathbf{R}_{0}^{T}=\mathbf{V}_{n}$ to obtain $\hat{\mathbf{A}}^{(n)}$.
\State  Normalize columns of $\hat{\mathbf{A}}^{(n)}$ to obtain factor matrix $\mathbf{A}^{(n)}$.
\State Recompute QR-decomposition for updated factor matrix $\mathbf{A}^{(n)}=\mathbf{Q}_{n}\mathbf{R}_{n}$.     
\EndFor
\State Compute the $fitness$.
\If{$fitness\ge \epsilon$}
\State Break
\EndIf
\EndFor
\end{algorithmic}
{\bfseries Output:} Matrices $\{\mathbf{A}^{(1)}, \cdots, \mathbf{A}^{(N)}\}$.
\label{ALS-QR-BRE}
\end{algorithm}

At the initial stage of the Algorithm \ref{ALS-QR-BRE}, if the discrepancy between the fitted values in the first two iterations is large, the extrapolation operation for $\mathbf{Q}_0$ in line 5 is temporarily omitted to prevent excessive fluctuations. During this phase, the Algorithm \ref{ALS-QR-BRE} reduces to ALS-QR-BR, performing iterative tensor computation without extrapolation. Once the variation between successive fitted values falls below a predefined threshold, the extrapolation is activated, accelerating convergence and improving both computational efficiency and fitting accuracy.

\subsection{Complexity analysis}\label{Cost}
The complexity analysis primarily focuses on optimizing multiple TTM operations through restructured dimension tree. To evaluate its effectiveness in reducing computational complexity, we compare three computational scenarios over the first three iterations:
\begin{itemize} 
\item Without the dimension tree (matching to line 5 of Algorithm \ref{CP-ALS-QR}).
\item Standard dimension tree ( hierarchical decomposition). 
\item Restructured dimension tree (corresponding to line 6 of Algorithm \ref{ALS-QR-BRE}). 
\end{itemize} 

We derive theoretical complexity bounds for these three cases, focusing on the cumulative computational cost of line 6 in Algorithm \ref{ALS-QR-BRE} during the first three iterations. The complexity of other steps in Algorithm \ref{ALS-QR-BRE} has been analyzed in \cite{RM20} and is not repeated here. Our study takes third-order and fourth-order tensors as examples and subsequently extends the computational complexity analysis to $N$th-order tensors.

To clarify the computational complexity of sequential TTM operations, we consider a third-order tensor $\mathcal{X} \in \mathbb{R}^{I_1 \times I_2 \times I_3}$  with the CP decomposition rank $R$. The first TTM along mode-$1$ incurs a computational cost of $2I_1I_2I_3R$, followed by the second TTM along mode-$2$ with a cost of  $2I_2I_3R^2$, and the third TTM along mode-$3$ requires $2I_3R^3$ operations. By extending this pattern, we derive the subsequent computational complexity results presented in the figures and tables. Table  \ref{tab_1} illustrates the dimension tree, while Table \ref{table_2} depicts the restructured dimension tree.

\begin{table}[!ht]
\begin{minipage}[p]{0.5\textwidth} 
\centering 
\includegraphics[width=4cm,height =4cm]{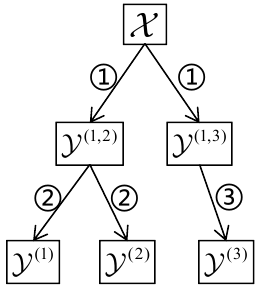} 
\end{minipage}
\begin{minipage}[p]{0.5\textwidth}
\centering
\begin{tabular}{|c|c|}
\hline		
Number   & Cost \\ \hline
\textcircled{1} & $2I_1I_2I_3R \times 2$
\\ \hline
\textcircled{2} &$2I_1I_2R^2 \times 2$ \\ \hline
\textcircled{3}& $2I_1I_3R^2 $\\ \hline
\end{tabular}
\label{tab:my_label}
\end{minipage}
\caption{The left diagram depicts the iterative structure of the dimension tree for third-order tensor, while the right table summarizes the total computational complexity per iteration. Identical-index edges indicate equivalent computational complexity.}
\label{tab_1}
\end{table}

\begin{table*}[!ht]
   \centering 
   \includegraphics[width=14cm,height = 4cm]{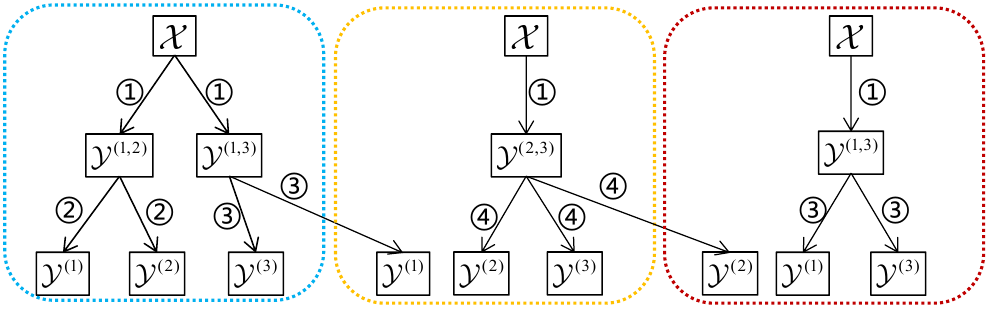}
   \label{fig_5}
   \vspace{10pt}
   \begin{tabular*}{\textwidth}{@{\extracolsep\fill}ccccc@{\extracolsep\fill}}
   \toprule
   \textbf{Number} & \textbf{\textcircled{1}}  & \textbf{\textcircled{2}}  & \textbf{\textcircled{3}}& \textbf{\textcircled{4}} \\
   \midrule
   Cost & $2I_1I_2I_3R \times 4$  & $2I_1I_2R^2 \times 2$ & $2I_1I_3R^2 \times 4$ &  $2I_2I_3R^2 \times 3$\\
  \bottomrule
  \end{tabular*}
  \vspace{5pt}
   \caption{
   The diagram above illustrates the iterative structure of branch reutilization of dimension tree for third-order tensor, while the table below shows the total computational complexity for each index throughout the iteration.}
  \label{table_2}
\end{table*}

Through synthesizing the iterative complexities of the dimension tree (Table \ref{tab_1}), the restructured dimension tree (Table \ref{table_2}), and the theoretical foundations of TTM complexity, we analytically derive the total computational complexity for the first three iterations of a third-order tensor across the three scenarios presented in Table \ref{tab_3}.

\begin{table*}[!ht]
\centering
\setlength{\abovecaptionskip}{0pt}
\setlength{\belowcaptionskip}{0pt}
\caption{Computational complexity of three scenarios of third-order tensor.\label{tab_3}}
\fontsize{7}{7}\selectfont
\begin{tabular*}{\textwidth}{@{\extracolsep\fill}cccc@{\extracolsep\fill}}
\toprule
\textbf{Scenario} & \textbf{Without the dimension tree}  & \textbf{Dimension tree}  & \textbf{Restructured dimension tree} \\
\midrule
Cost & $18I_1I_2I_3R$  & $12I_1I_2I_3R$& $8I_1I_2I_3R$ \\ &$+6I_1I_2R^2+6I_2I_3R^2+6I_1I_3R^2$& $+12I_1I_2R^2+6I_1I_3R^2$ & $+4I_1I_2R^2+8I_1I_3R^2+6I_2I_3R^2$   \\
\bottomrule
\end{tabular*}
\end{table*}

When $I_1 = I_2 = I_3$, the term $I_1I_2I_3R$ dominates the overall computational cost. From the quantity of this term, we can conclude that the branch reutilization of  dimension tree significantly reduces computational complexity compared to the scenario without the dimension tree, while also outperforming the dimension tree. 

Let $\mathcal{X} \in \mathbb{R} ^{I_1 \times I_2 \times I_3 \times I_4}$ be a fourth-order tensor, and assume the CP approximation rank of $R$. By leveraging the theoretical framework TTM complexity, we derive the results of the dimension tree presented in Table \ref{tab_4}, while Table \ref{tab_5} presents the results of the restructured dimension tree.

\begin{table}[!ht]
\begin{minipage}[p]{0.5\textwidth} 
\centering 
\includegraphics[scale=0.8]{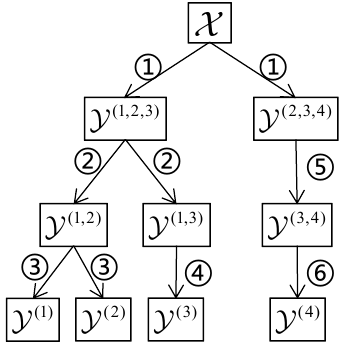} 
\end{minipage}
\begin{minipage}[p]{0.5\textwidth}
\centering
\begin{tabular}{|c|c|}
\hline		
Number   & Cost \\ \hline
\textcircled{1} & $2I_1I_2I_3I_4R \times 2$
\\ \hline
\textcircled{2} &$2I_1I_2I_3R^2 \times 2$ \\ \hline
\textcircled{3}& $2I_1I_2R^3 \times 2 $\\ \hline
\textcircled{4} &$2I_1I_3R^3 $   \\ \hline
\textcircled{5} &$2I_2I_3I_4R^2$   \\ \hline
\textcircled{6} & $2I_3I_4R^3$  \\ \hline
\end{tabular}
\end{minipage}
\caption{The left diagram illustrates the iterative structure of the dimension tree for a fourth-order tensor, while the table on the right summarizes the corresponding computational complexity of each index in the dimension tree.}
\label{tab_4}
\end{table}

\begin{table}[!ht]
\centering 
\includegraphics[width=14cm,height = 4.8cm]{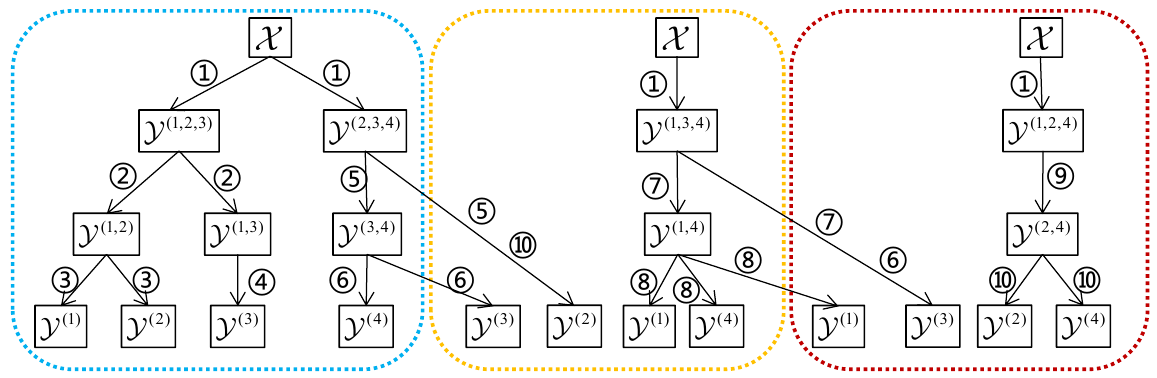}
\label{10pt}
\vspace{10pt}
\centering
\begin{tabular*}{\textwidth}{@{\extracolsep\fill}ccccccc@{\extracolsep\fill}}
\toprule
\textbf{Number} & \textbf{\textcircled{1}}  & \textbf{\textcircled{2}}  & \textbf{\textcircled{3}}& \textbf{\textcircled{4}} &\textbf{\textcircled{5}}\\
\midrule
Cost & $2I_1I_2I_3I_4R \times 4$  & $2I_1I_2I_3R^2 \times 2$ & $2I_1I_2R^3\times 2$  &  $2I_1I_3R^3$& $2I_2I_3I_4R^2\times2$ \\ \hline
\toprule
\textbf{Number} & \textbf{\textcircled{6}}  & \textbf{\textcircled{7}}  & \textbf{\textcircled{8}}& \textbf{\textcircled{9}} &\textbf{\textcircled{10}}\\ \hline
Cost &$2I_3I_4R^3\times3$ &$2I_1I_3I_4R^2 \times 2$ &$2I_1I_4R^3\times3$ &$2I_1I_2I_4R^2$ &$2I_2I_4R^3\times3$ \\
\bottomrule
\end{tabular*}
\caption{The figure above shows the iterative structure of restructured dimension tree for the fourth-order tensor, while the table below presents the total computational complexity for each index during the iteration.}\label{tab_5}
\end{table}

The iterative complexity of the dimension tree (Table \ref{tab_4}), the complexity of the restructured dimension tree ( Table \ref{tab_5}), and the TTM complexity framework are integrated to enable the systematic deduction of the total computational complexity for the three scenarios, with the results presented in Table \ref{tab_6}.

\begin{table*}[!ht]
\centering
\setlength{\abovecaptionskip}{0pt}
\setlength{\belowcaptionskip}{0pt}
\caption{Computational complexity of three scenarios of fourth-order tensor.\label{tab_6}}
\fontsize{7}{7}\selectfont
\begin{tabular*}{\textwidth}{@{\extracolsep\fill}cccc@{\extracolsep\fill}}
\toprule
\textbf{Scenario} & \textbf{Without the dimension tree}  & \textbf{Dimension tree} & \textbf{Restructured dimension tree} \\
\midrule
Cost & $24I_1I_2I_3I_4R+12I_1I_2I_3R^2$  & $12I_1I_2I_3I_4R+12I_1I_2I_3R^2$ & $8I_1I_2I_3I_4R+4I_2I_4R^3++6I_1I_4R^3$\\ 
 &$+12I_2I_3I_4R^2+12I_1I_2R^3$  & $+6I_2I_3I_4R^2+12I_1I_2R^3$ & $+2I_1I_3R^3+6I_3I_4R^3+4I_1I_2R^3$   \\
 & $+12I_3I_4R^3$& $+6I_1I_3R^3+6I_3I_4R^3$ & $+2I_2I_3R^3+4I_1I_3I_4R^2+2I_1I_2I_4R^2$ \\ & & &$+4I_2I_3I_4R^2+4I_1I_2I_3R^2$ \\
\bottomrule
\end{tabular*}
\vspace{-10pt}
\end{table*}

When $I_1=I_2=I_3=I_4$, the dominant computational cost is $I_1I_2I_3I_4R$. In restructured dimension tree, the occurrence of this term is fewer compared to both the classical dimension tree and without the dimension tree. Consequently, the computational cost of the restructured dimension tree is lower than that of the other two scenarios.

\begin{remark}
The reconstructed dimension tree reduces the main computational burden. Meanwhile, by distributing the remaining computational load in a non-centralized manner across various indices, it enhances the algorithm's adaptability to different datasets.
\end{remark}

We extend the analysis to 
$N$-th order tensors and, under the assumption of  
$I=I_1=I_2=\cdots=I_N$, obtain the computational cost of TTM for $N$-th order tensors in the first three iterations of the algorithm. The result is presented in Corollary \ref{cor_1}. 

\vspace{8pt}
\begin{corollary}\label{cor_1}
\raggedright
Computational complexity of three scenarios of $N$th-order tensor $(N>4)$.
\fontsize{7}{7}\selectfont
\begin{tabular*}{\textwidth}{@{\extracolsep\fill}cccc@{\extracolsep\fill}}
\toprule
\textbf{Scenario} & \textbf{Without the dimension tree}  & \textbf{Dimension tree}  & \textbf{Restructured dimension tree} \\
\midrule
Cost & $6NI^{N}R$ & $12I^{N}R$ & $8I^{N}R$\\  
     & $+6NI^{N-1}R^2+\cdots+6NI^{2}R^{N-1}$ & $+18I^{N-1}R^2+\cdots+6NI^{2}R^{N-1}$  & $+14I^{N-1}R^2+\cdots+6NI^{2}R^{N-1}$\\
\bottomrule
\end{tabular*}
\end{corollary}

The above analysis provides the total computational cost for $N$-order tensors under three distinct scenarios: without the dimension tree,  dimension tree, and restructured dimension tree. A comparative examination of the dominant computational terms reveals that the restructured dimension tree yields a 33\% reduction in computational cost and demonstrates superior adaptability to diverse data structures compared with the dimension tree.

\section{Numerical experiments}\label{numerical} 
In this section, we evaluate the performance of our proposed algorithms, ALS-QR-BRE and ALS-QR-BR, on random, synthetic, and real tensors, and compare them with the classical CP-ALS \cite{RM23}, CP-ALS-PINV, the unmodified CP-ALS-QR \cite{RM20}, CP-ALS-QR-SVD and CP-GEVD \cite{Tl} (The generalization of the eigenvalue method for CP decomposition \cite{leurgans1993threeway,Sanchez1990}) algorithms to assess their efficiency. All numerical experiments were conducted on a PC with an i7-12700 processor, 16 GB of RAM, and implemented in MATLAB using Tensortoolbox \cite{A4} and Tensorlab \cite{Tl}. A comparison of each algorithm is provided in the table \ref{tab_8}. This implementation will be made available at \url{https://github.com/youmengx658/ALS-QR-BRE} upon acceptance of this manuscript.
\begin{table}[!ht]
    \centering 
    \setlength{\abovecaptionskip}{0pt}
    \setlength{\belowcaptionskip}{0pt}
    \caption{Algorithm.}
    \begin{tabular}{c|c}
         \hline		
          Name   & Interpretation \\ \hline
       CP-ALS (ALS)& Algorithm \ref{CP-ALS} (line 6) with Cholesky for the inverse. 
      \\ \hline
       CP-ALS-PINV (PINV)& Algorithm \ref{CP-ALS} (line 6) with SVD  for the inverse. \\ \hline
        CP-GEVD  &    The generalized eigenvalue method for CP decomposition.     \\ \hline
      CP-ALS-QR (QR)& Algorithm \ref{CP-ALS-QR}.\\ \hline
      CP-ALS-QR-SVD (QR-SVD)& Algorithm \ref{CP-ALS-QR} (line 7) with SVD for the inverse.\\  \hline
       CP-ALS-QR-DT (QR-DT)& Algorithm \ref{CP-ALS-QR} (line 6) with the dimension tree.\\  \hline
        ALS-QR-BR(QR-BR)& Algorithm \ref{ALS-QR-BRE} with $\beta=0$.\\  \hline
         ALS-QR-BRE (QR-BRE)& Algorithm \ref{ALS-QR-BRE}.\\  \hline
    \end{tabular}
    \label{tab_8}
\end{table}

The abbreviations in parentheses in Table \ref{tab_8} represent the algorithms used for the evaluation of the average iteration time in Subsection \ref{test_random}. The full names of the algorithms are employed in the experiments involving synthetic and real tensors, as discussed in Subsections \ref{sy_tensor} and \ref{real_tensor}.
\vspace{-5pt}
\subsection{Approximation error}\label{error}
The following experiments are conducted based on the approximation error defined as follows:
$$
\sqrt{\left \| \mathcal{X} -\mathcal{K}  \right \|^2 } =\sqrt{\left \|\mathcal{X}   \right \|^2 -2\left \langle \mathcal{X}, \mathcal{K}\right \rangle +\left \|\mathcal{K}  \right \|^2 },
$$
where  $\left \| \mathcal{X}  \right \|^2 $  is precomputed and remains constant throughout the iteration process. $\left \langle \mathcal{X} ,\mathcal{K}  \right \rangle $ and $\left \| \mathcal{K}  \right \|^2$ are computed using intermediate quantities generated during algorithm iterations.

In the CP-ALS algorithm, the approximation $\mathcal{K}$ is expressed as $\mathbf{K}_{(n)}=\hat{\mathbf{A}}^{(n)} \mathbf{P}^{(n)T}$, where $\hat{\mathbf{A}}^{(n)}=\mathbf{A}^{(n)} \cdot diag (\lambda)$ and $\mathbf{P}^{(n)}=\mathbf{A}^{(1)} \odot \cdots \odot \mathbf{A}^{(n-1)} \odot \mathbf{A}^{(n+1)} \odot  \cdots \odot \mathbf{A}^{(N)}$. Consequently, the inner product between the tensor $\mathcal{X}$ and the approximation $\mathcal{K}$ is given by:
$$
\left \langle \mathcal{X},\mathcal{K}  \right \rangle =\left \langle \mathbf{X}_{(n)},\mathbf{K}_{(n)} \right \rangle=\left \langle \mathbf{X}_{(n)},\hat{\mathbf{A}}^{(n)} \mathbf{P}^{(n)T} \right \rangle=\left \langle \mathbf{X}_{(n)} \mathbf{P}^{(n)},\hat{\mathbf{A}}^{(n)} \right \rangle =\left \langle \mathbf{M}_{n} ,\hat{\mathbf{A}}^{(n)} \right \rangle, 
$$
where $\mathbf{M}_{n}$ computed as described in line 5 of Algorithm \ref{CP-ALS}. Similarly, the squared norm of $\mathcal{K}$ is calculated as:
$$
\left \| \mathcal{K}  \right \|^2= \left \langle \hat{\mathbf{A}}^{(n)} \mathbf{P}^{(n)T},\hat{\mathbf{A}}^{(n)} \mathbf{P}^{(n)T} \right \rangle=\left \langle \mathbf{P}^{(n)T}\mathbf{P}^{(n)},(\hat{\mathbf{A}}^{(n)})^T\hat{\mathbf{A}}^{(n)}  \right \rangle =\left \langle \Gamma ^{(n)},diag(\mathbf{\lambda}  )\mathbf{S}_n diag(\mathbf{\lambda} ) \right \rangle ,
$$
where $\Gamma ^{(n)}$ is computed in line 3 of Algorithm \ref{CP-ALS} and the $\mathbf{S}_n=(\mathbf{A}^{(n)})^T\mathbf{A}^{(n)}$ is gram matrix of the $n$th factor matrix. Typically, the index $n$ refers to the last mode in the subiteration. This approach for calculating the approximation error is well-established in  \cite{A4,Ballard2025TensorDecompositions, A1, A2,Tb,A3}. 

According to the previously described computation method, $\mathbf{K}_{(n)}=\hat{\mathbf{A}}^{(n)} \mathbf{P}^{(n)T}$ where $\mathbf{P}^{(n)}=(\mathbf{Q}_N\otimes \cdots \otimes \mathbf{Q}_{n+1}\otimes \mathbf{Q}_{n-1}\otimes \cdots \otimes \mathbf{Q}_{1})\mathbf{Q}_0\mathbf{R}_0$ for the Algorithm \ref{CP-ALS-QR}. Consequently, the inner product between $\mathcal{X}$ and 
$\mathcal{K}$ is expressed as:
$$
\left \langle \mathcal{X},\mathcal{K}  \right \rangle =\left \langle \mathbf{X}_{(n)},\mathbf{K}_{(n)} \right \rangle=\left \langle \mathbf{X}_{(n)},\hat{\mathbf{A}}^{(n)} \mathbf{P}^{(n)T} \right \rangle=\left \langle \mathbf{X}_{(n)} \mathbf{P}^{(n)},\hat{\mathbf{A}}^{(n)} \right \rangle =\left \langle \mathbf{V}_{n} ,\hat{\mathbf{A}}^{(n)} \mathbf{R}_0^T \right \rangle,
$$
where $\mathbf{V}_{n}$ is computed as outlined in line 6 of Algorithm \ref{CP-ALS-QR}, and $\mathbf{R}_0$ 
is the upper triangular matrix obtained from the QR decomposition in line 4 of Algorithm \ref{CP-ALS-QR}. Similarly, the squared norm of $\mathcal{K}$ is given by:
$$
\left \| \mathcal{K}  \right \|^2= \left \langle \hat{\mathbf{A}}^{(n)} \mathbf{P}^{(n)T},\hat{\mathbf{A}}^{(n)} \mathbf{P}^{(n)T} \right \rangle=\left \langle \mathbf{P}^{(n)T}\mathbf{P}^{(n)},(\hat{\mathbf{A}}^{(n)})^T\hat{\mathbf{A}}^{(n)}  \right \rangle =\left \langle \mathbf{R}_0^T \mathbf{R}_0,diag(\lambda )\mathbf{R}_n^T \mathbf{R}_ndiag(\lambda ) \right \rangle ,
$$
    where $\mathbf{R}_n$ is the upper triangular matrix in the QR decomposition of $\mathbf{A}^{(n)}$. This error calculation approach for the CP-ALS-QR algorithm was initially introduced in \cite{RM20}. The fitness metric, which evaluates the quality of the approximation, is defined as:
$$
Fitness = 1-\frac{\sqrt{\left \| \mathcal{X} -\mathcal{K}  \right \|^2 }}{\left \| \mathcal{X} \right \| },
$$
where a fitness value approaching 1 indicates a better approximation of $\mathcal{X}$ by $\mathcal{K}$.

\subsection{Evaluate the average iteration time on random tensors}\label{test_random}
\begin{figure}[!ht]
\vspace{-5pt}
 \begin{minipage}[p]{0.5\textwidth}
\centering  
\subfloat{
\label{Fig.sub.01}
\includegraphics[scale=0.38]{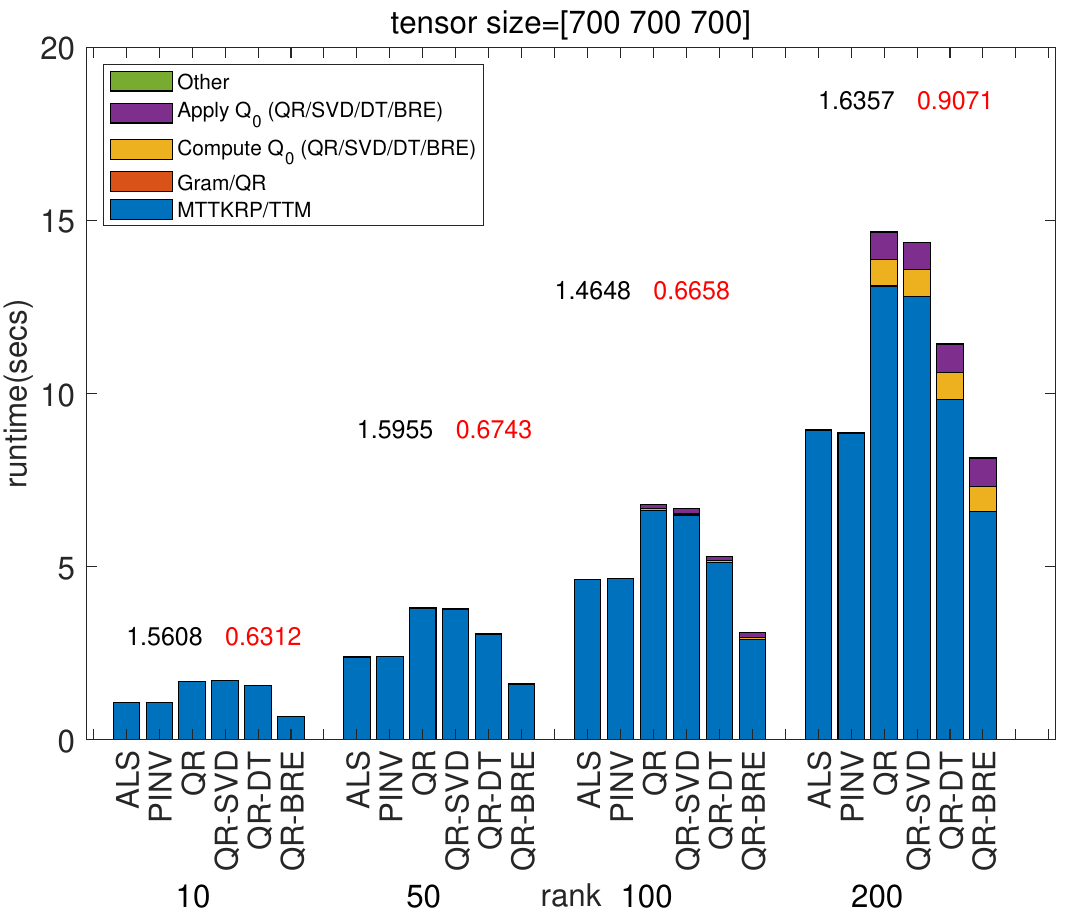}}
\vspace{5pt}
\end{minipage}
\begin{minipage}[p]{0.5\textwidth}
\centering
\subfloat{
\label{Fig.sub.02}
\includegraphics[scale=0.38]{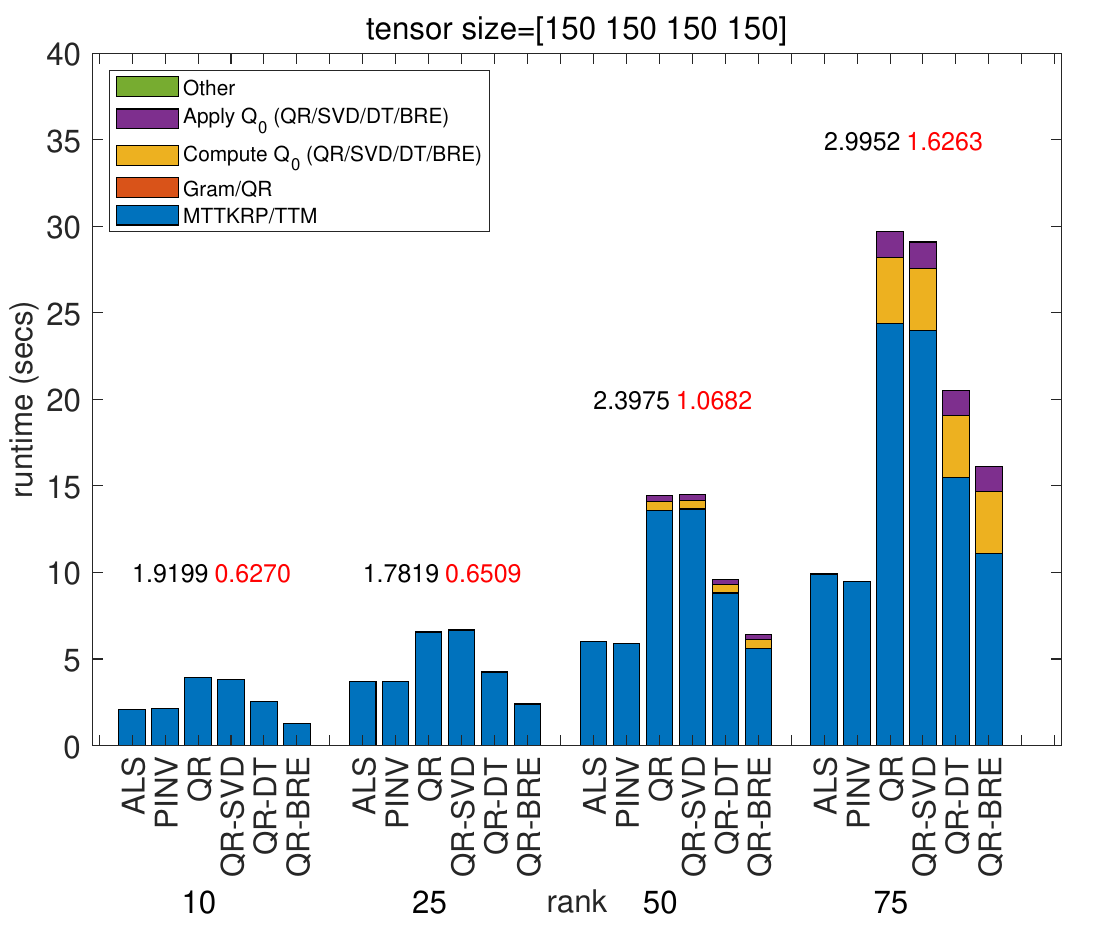}}
\vspace{5pt}
\end{minipage}
\begin{minipage}{0.5\textwidth}
\centering
\subfloat{
\label{Fig.sub.3}
\includegraphics[scale=0.38]{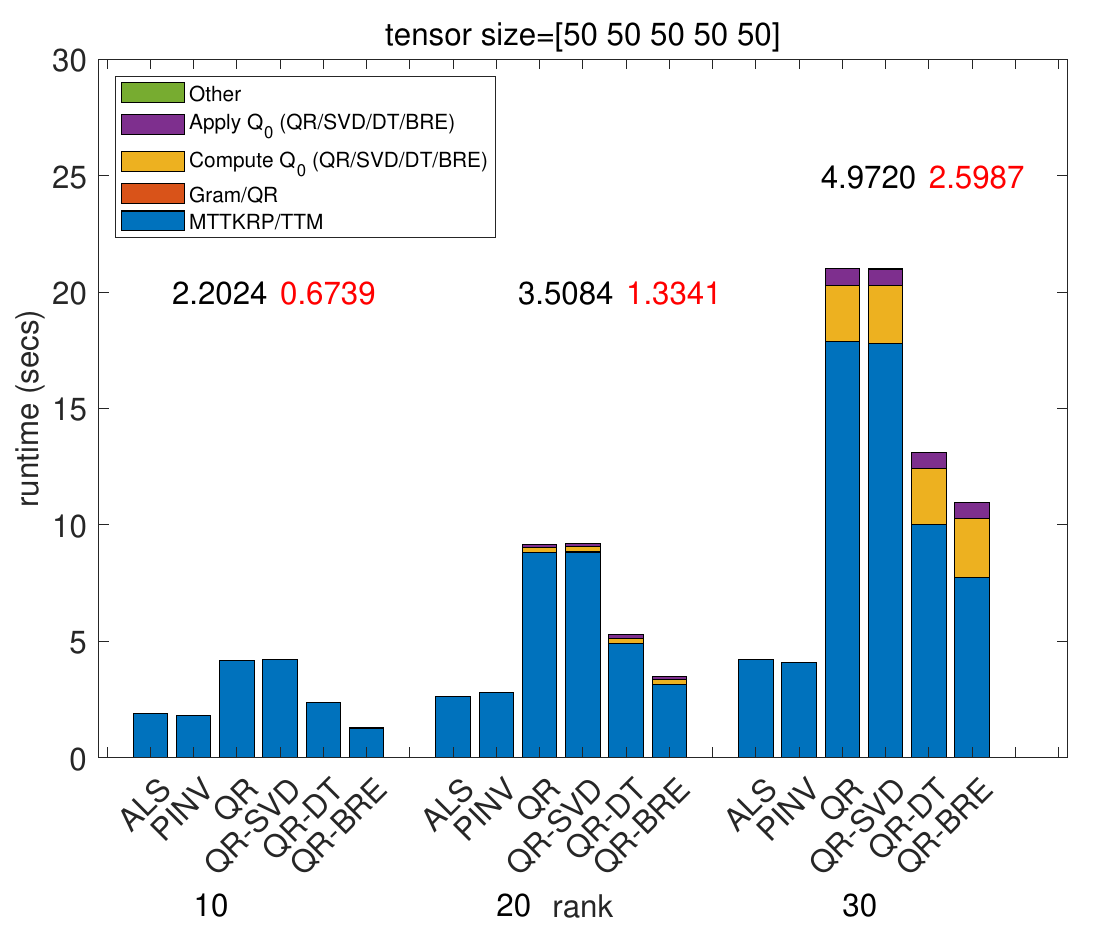}}
\end{minipage}
\begin{minipage}{0.5\textwidth}
\centering
\subfloat{
\label{Fig.sub.4}
\includegraphics[scale=0.38]{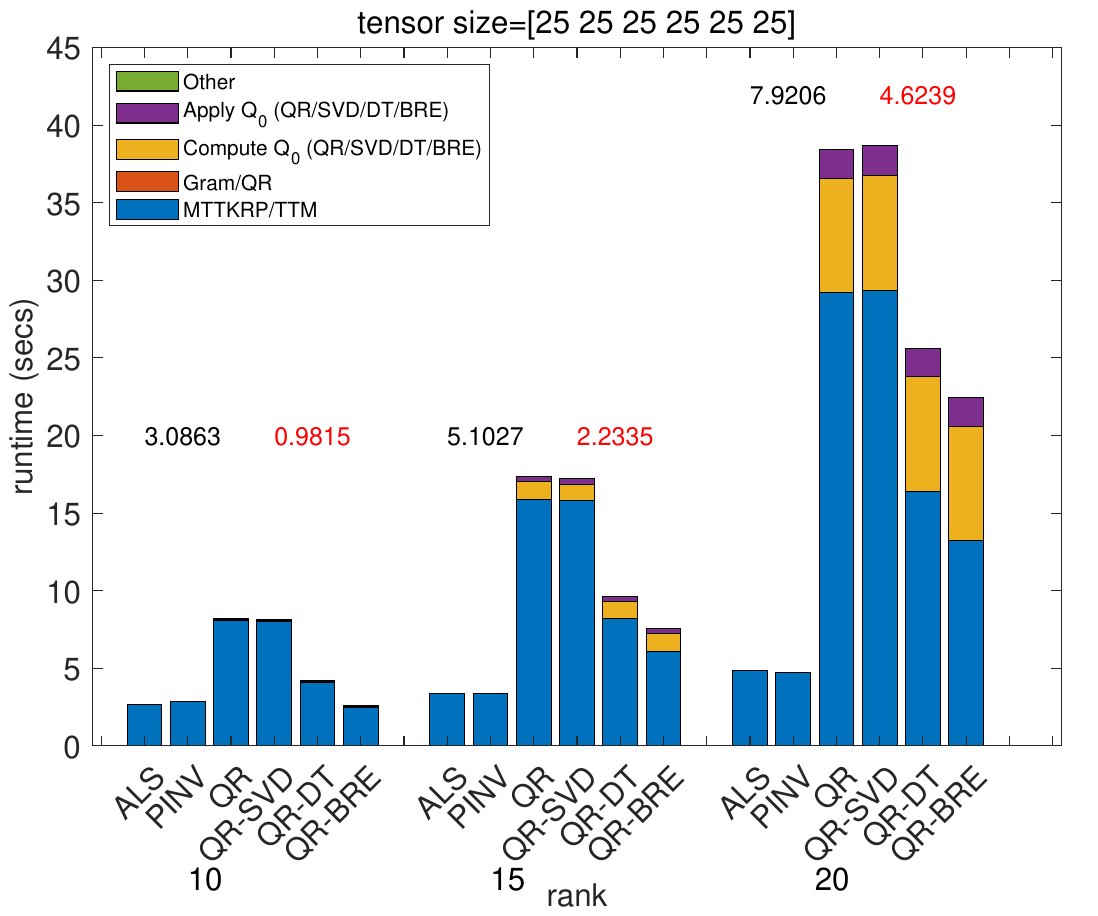}}
\end{minipage}
\caption{The average runtime (in seconds) for a single iteration of ALS, PINV, QR, QR-SVD, QR-DT, and QR-BRE is presented as the rank increases, for a third-order tensor of size 700 (top left), a fourth-order tensor of size 150 (top right), a five-order tensor of size 50 (bottom left), and a six-order tensor of size 25 (bottom right).}
\label{fig-7}
\end{figure}

In this subsection, we conduct experiments on random tensors of orders three, four, five, and six to evaluate the impact of increasing rank on the average iteration time of each algorithm (excluding the CP-GEVD algorithm). Specifically, the dimensions of the third-order tensor are set to 700 for each mode, while those of the fourth-order, fifth-order, and sixth-order tensors are set to 150, 50, and 25, respectively. Each algorithm is executed 10 times. To minimize the influence of initialization bias, the runtime of the first execution is discarded, and the average runtime of the subsequent executions is computed. The experimental results for tensors of varying ranks across these four cases are presented in Figure \ref{fig-7}.

As shown in Figure \ref{fig-7}, MTTKRP and Gram (refer to lines 3 and 5 in Algorithm \ref{CP-ALS}) represent the proportion of computational time spent on calculating $\mathbf{M}_{i}=\mathbf{X}_{(i)}\mathbf{P}^{(i)}$ and computing the Gram matrix for each factor matrix, respectively, within the ALS and PINV algorithms. TTM and QR denote the computational time consumed by TTM operations and QR decomposition in the QR, QR-SVD, QR-DT, and QR-BRE algorithms. Specifically, the TTM and QR operations correspond to lines 5 and initialize in Algorithm \ref{CP-ALS-QR}, as well as lines 6 and initialize in Algorithm \ref{ALS-QR-BRE}. The computation of $\mathbf{Q}_0$ involves performing QR decomposition on the Khatri-Rao product $\mathbf{Z}_n$ (refer to line 4 in Algorithms \ref{CP-ALS-QR} and \ref{ALS-QR-BRE}), while the application of $\mathbf{Q}_0$ corresponds to a matrix multiplication operation (refer to line 6 in Algorithm \ref{CP-ALS-QR} and line 7 in Algorithm \ref{ALS-QR-BRE}). Other steps, such as extrapolation, weight computation, and error calculation, are categorized as "other" due to their minimal contribution to the overall computational time. Since the time spent on extrapolation is negligible, the iteration time for the QR-BR algorithm is the same as that for QR-BRE. Additionally, for each rank, the runtime ratio between ALS and QR-BRE is highlighted in red, while the runtime ratio between ALS and QR is highlighted in black.

From the experimental results in Figure \ref{fig-7}, it can be observed that the primary computational cost per iteration is MTTKRP for the ALS and PINV algorithms, whereas for the QR, QR-SVD, QR-DT, and QR-BRE algorithms, the dominant cost is TTM. For third-order tensors, the proposed QR-BRE algorithm consistently outperforms the other five algorithms in iteration speed across different ranks. The time proportion of TTM indicates that the restructured dimension tree provides superior acceleration compared to the dimension tree. The runtime ratio in the red line further confirms that the proposed QR-BRE algorithm consistently achieves shorter runtimes than the ALS algorithm. For fourth-order, fifth-order, and sixth-order tensors, the QR-BRE algorithm accelerates the QR algorithm regardless of rank. When the rank is small, it also outperforms the ALS algorithm in iteration speed, as reflected in the red runtime ratio in Figure \ref{fig-7}. Its significant computational advantage over the dimension tree is evident from the TTM proportion in the QR-DT and QR-BRE algorithms.

\subsection{Synthetic tensors}\label{sy_tensor}

This subsection primarily investigates synthetic tensors with a consistent true rank while varying factor collinearity, tensor dimensions, and noise levels. Specifically, two different sizes of third-order and fourth-order tensors are generated, each with a true rank of $\hat{R}$, different levels of factor collinearity, and varying noise conditions. The collinearity in each mode is represented by $c_i\in(0,1)$. Homoscedastic noise is introduced with an intensity controlled by the parameter $l_1$. When $l_2 > 0$, heteroscedastic noise is additionally introduced into the tensor, with $l_2$ controlling its proportion.


These tensors are generated following the procedures outlined in references \cite{A6, evans2023blockwise}, and \cite{tomasi2006comparison}. The main steps are as follows: First, specify the true rank $\hat{R}$, tensor size, parameters $l_1$ and $l_2$, and the collinearity coefficients $c_i$. Second, generate an $R \times R$ matrix $\mathbf{K}$, where the diagonal elements are all 1, and the off-diagonal elements are the collinearity coefficient $c_i$. Third, compute the Cholesky factor $\mathbf{C}$ of the $\mathbf{K}$. Fourth, create a matrix $\mathbf{M}$ of size $I \times R$ with elements sampled from a standard normal distribution, and perform column-wise orthogonal normalization to obtain $\mathbf{Q}$. Finally, generate the $n$th factor matrix $\hat{\mathbf{B}}^{(n)} = \mathbf{Q}\mathbf{C}$. This process generates multiple factor matrices $\hat{\mathbf{B}}^{(n)}$, which are then used to construct the tensor $\mathcal{X}$. For $0 \leq l_1 < 100$, homoscedastic noise is added to the tensor $\mathcal{X}$ as follows.
$$
\mathcal{X}' =\mathcal{X} +\frac{1}{\sqrt{100/l_1-1}} \frac{\left \|\mathcal{X}   \right \|_F }{\left \|  \mathcal{N}_1 \right \|_F } \mathcal{N}_1.
$$
The size of tensor $\mathcal{N}_1$ matches that of tensor $\mathcal{X}$, and the elements of tensor $\mathcal{N}_1$ are drawn from a standard normal distribution. When $l_1 = 0$, $\mathcal{X}' = \mathcal{X}$, and when $l_1$ is close to 100, $\mathcal{X}'$ approximates random noise.When $l_2 > 0$, heteroscedastic noise is added to the tensor $\mathcal{X}'$ as follows.
$$
\mathcal{X}'' =\mathcal{X}' +\frac{1}{\sqrt{100/l_2-1}} \frac{\left \|\mathcal{X}'   \right \|_F }{\left \|  \mathcal{N}_2 \right \|_F } \mathcal{N}_2,
$$
where the elements of tensor $\mathcal{N}_2$ are drawn from a normal distribution with mean 0 and standard deviation 3. When $l_2=0$, $\mathcal{X}''=\mathcal{X}'$ and $l_2$ is approaches 100, $\mathcal{X}''$ is approximates random noise. The detailed parameters of the synthetic tensors used in the experiment can be found in Table \ref{tab-8}.

\begin{table*}[!ht]
\centering
\caption{The detailed parameters of the synthetic tensor.\label{tab-8}}
\begin{tabular*}{\textwidth}{@{\extracolsep\fill}ccccccc@{\extracolsep\fill}}
\toprule
\textbf{third-order} & \textbf{tensor size}  & \textbf{$c_1,c_2,c_3$}  & \textbf{ $\hat{R}$}& \textbf{$R$} &\textbf{$l_1$} &\textbf{$l_2$}\\ \hline
 1& [500 500 500] & [0.9 0.9 0.9] & 20&10 & 0.01 & 0\\ 
2& [500 500 500] & [0.9 0.9 0.9]  & 20 &75 & 0.01 & 0 \\ 
3& [500 500 500] & [0.9 0.9 0.9]  & 20 & 150& 0.01 &0 \\ 
4 & [600 600 600] & [0.5 .09 .09] & 20 & 10 &0.01 & 0.1 \\ 
 5& [600 600 600] & [0.5 .09 .09]& 20 & 75 &  0.01&0.1  \\ 
6 & [600 600 600] & [0.5 .09 .09] & 20 & 175& 0.01 & 0.1 \\ \hline
\toprule
\textbf{fourth-order} & \textbf{tensor size}  & \textbf{$c_1,c_2,c_3,c_4$}  & \textbf{ $\hat{R}$}& \textbf{$R$} &\textbf{$l_1$} &\textbf{$l_2$}\\ \hline
7& [100 100 100 100] &[0.9 0.9 0.9 0.9]& 20 & 10 & 0.1& 0 \\ 
8 & [100 100 100 100] &[0.9 0.9 0.9 0.9]& 20 &20  & 0.1 & 0 \\ 
9& [100 100 100 100] &[0.9 0.9 0.9 0.9]& 20& 30 &0.1  & 0 \\ 
 10 & [120 120 120 120]&[0.5 0.9 0.9 0.5]& 20 & 10 & 0.1 &0.01 \\ 
11& [120 120 120 120]&[0.5 0.9 0.9 0.5]& 20 & 20 & 0.1 & 0.01\\ 
12& [120 120 120 120]&[0.5 0.9 0.9 0.5]& 20 &  30& 0.1 & 0.01\\ 
\bottomrule
\end{tabular*}
\vspace{-10pt}
\end{table*}

In all experiments involving synthetic and real-world tensors, the two extrapolation hyperparameters for Algorithm \ref{ALS-QR-BRE} are selected according to the following rules: The hyperparameter $\alpha$ is fixed at $\frac{1}{10}$. The hyperparameter $\beta$ is chosen only if the difference in fitting performance between the first two iterations is less than 0.03. If the fitting coefficient at this point exceeds 0.95, $\beta$ is set to $\frac{1}{20000}$; if the fitting coefficient is between 0.90 and 0.95, $\beta$ is set to $\frac{1}{2000}$; if the fitting coefficient is between 0.70 and 0.90, $\beta$ is set to $\frac{1}{500}$; otherwise, $\beta$ is set to $\frac{1}{250}$. Once $\beta$ is selected, its value remains fixed throughout subsequent iterations.

\begin{figure}[!ht]
\vspace{-2pt}
\begin{minipage}[p]{0.5\textwidth}
\centering  
\subfloat[]{
\includegraphics[scale=0.38]{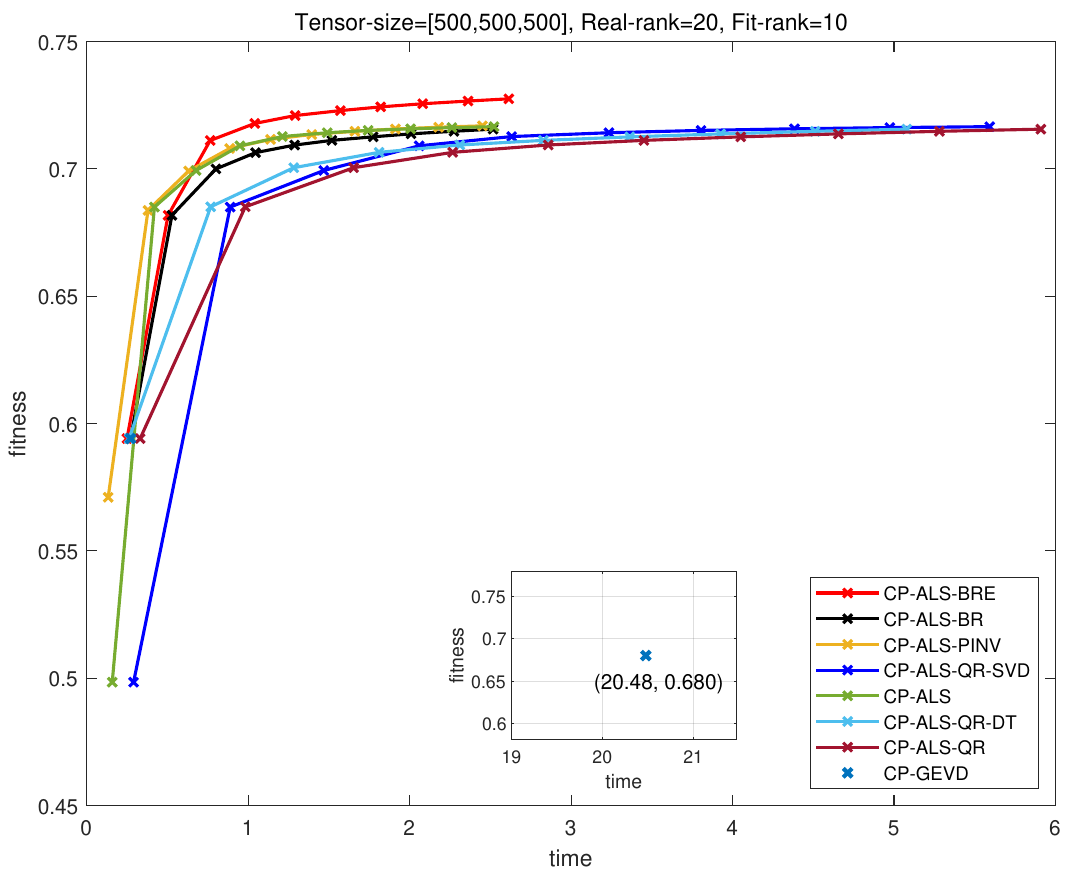}}
\end{minipage}
\begin{minipage}[p]{0.5\textwidth}
\centering
\subfloat[]{
\includegraphics[scale=0.38]{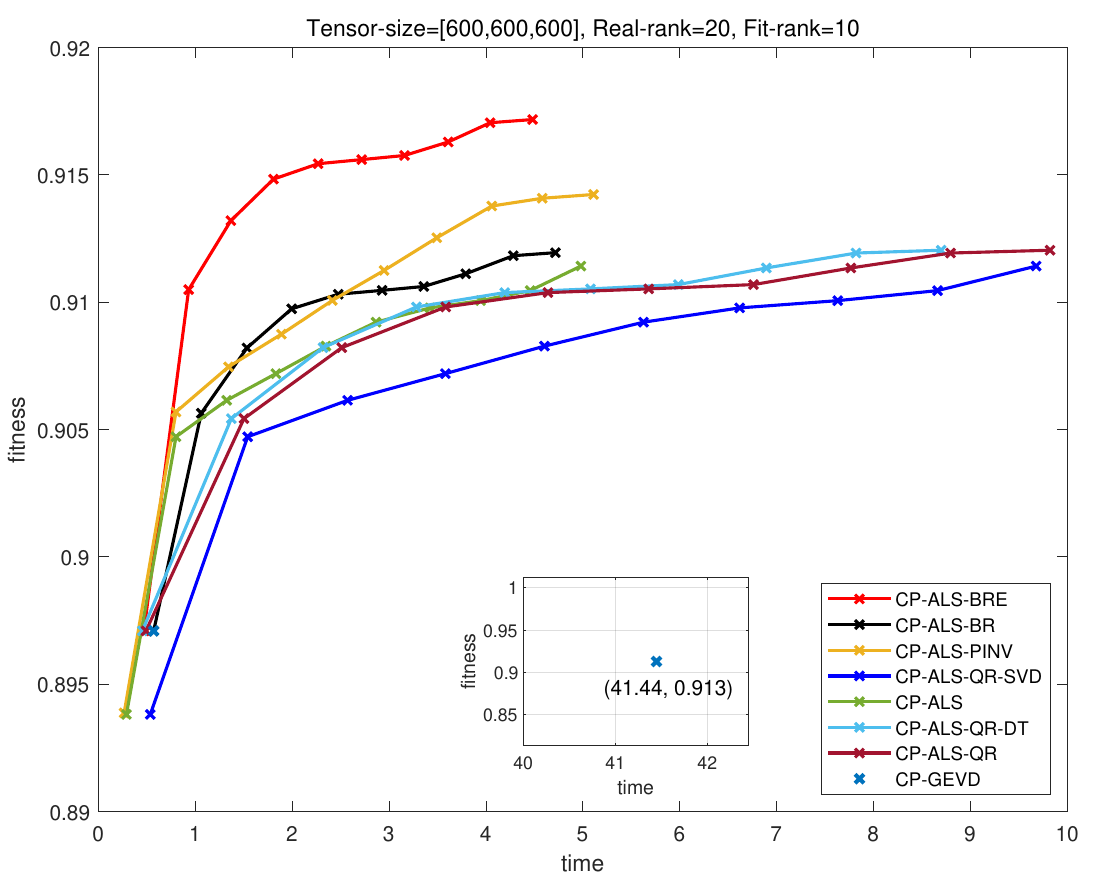}}
\end{minipage} 
\begin{minipage}{0.5\textwidth}
\centering
\subfloat[]{
\includegraphics[scale=0.38]{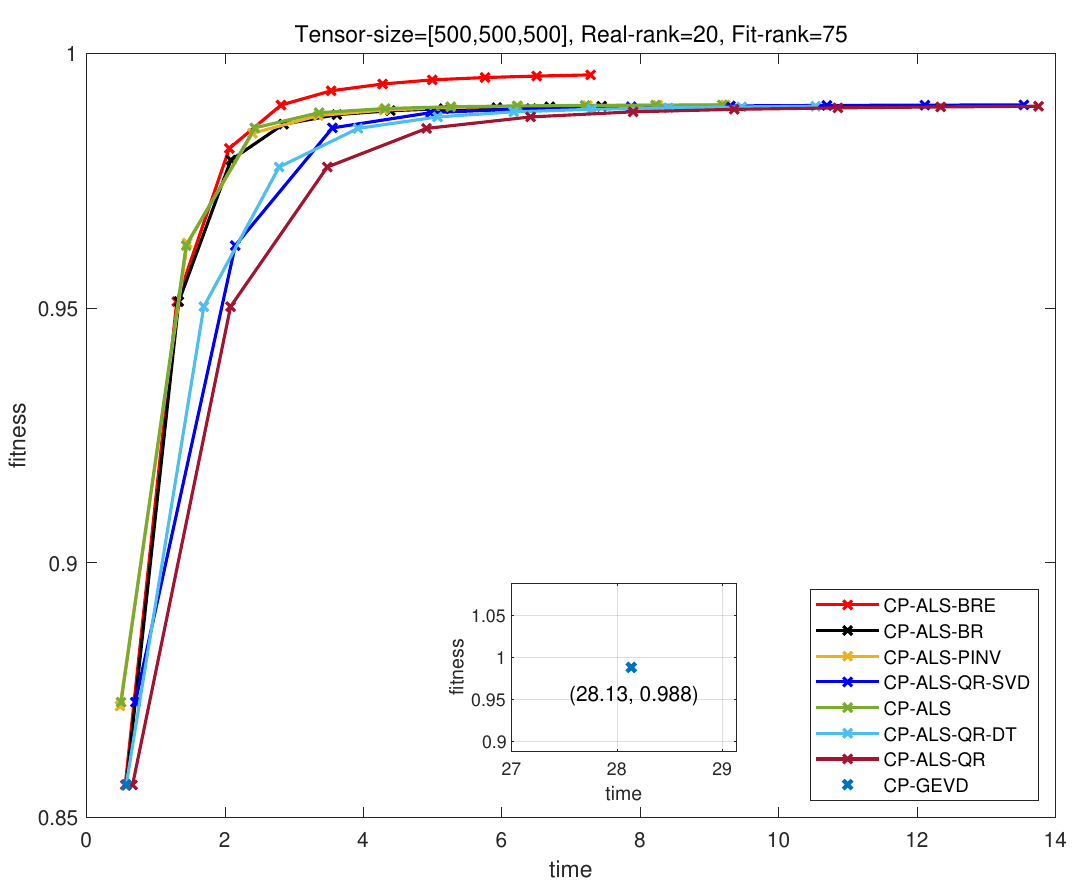}}
\end{minipage}
\begin{minipage}{0.5\textwidth}
\centering
\subfloat[]{
\includegraphics[scale=0.38]{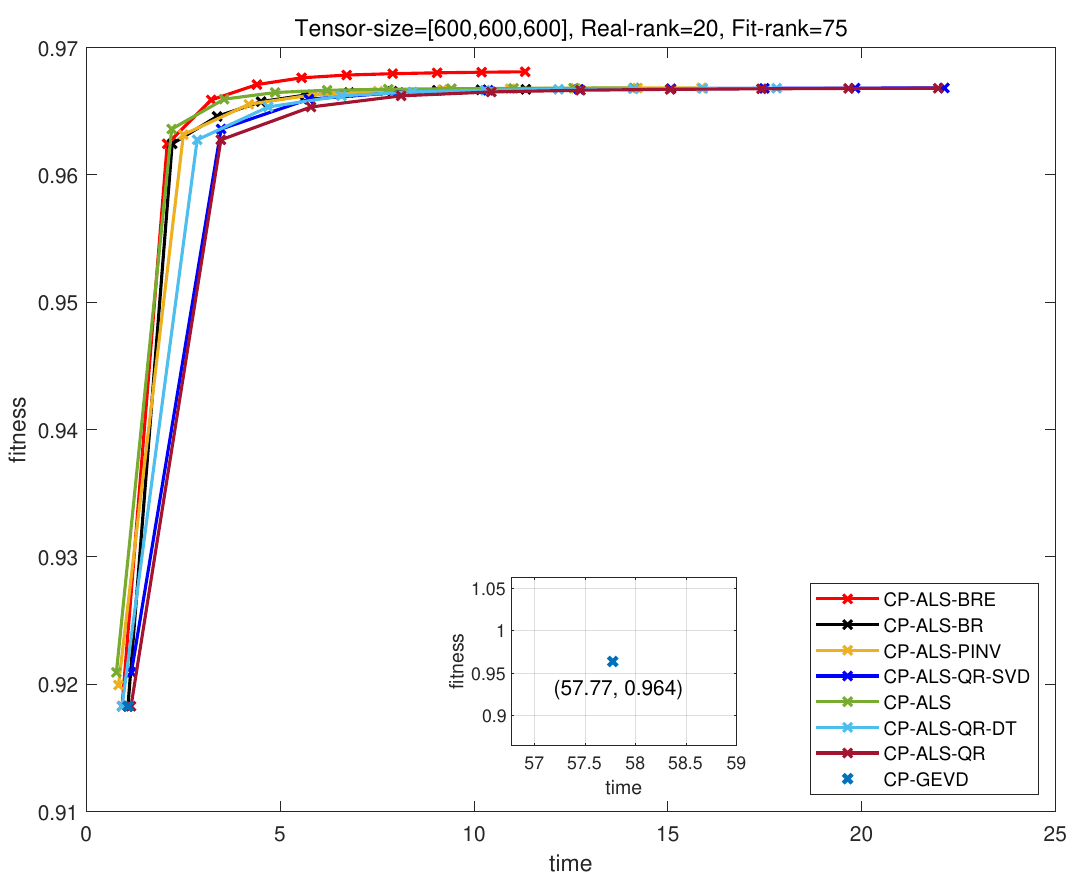}}
\end{minipage}
\begin{minipage}{0.5\textwidth}
\centering
\subfloat[]{
\includegraphics[scale=0.38]{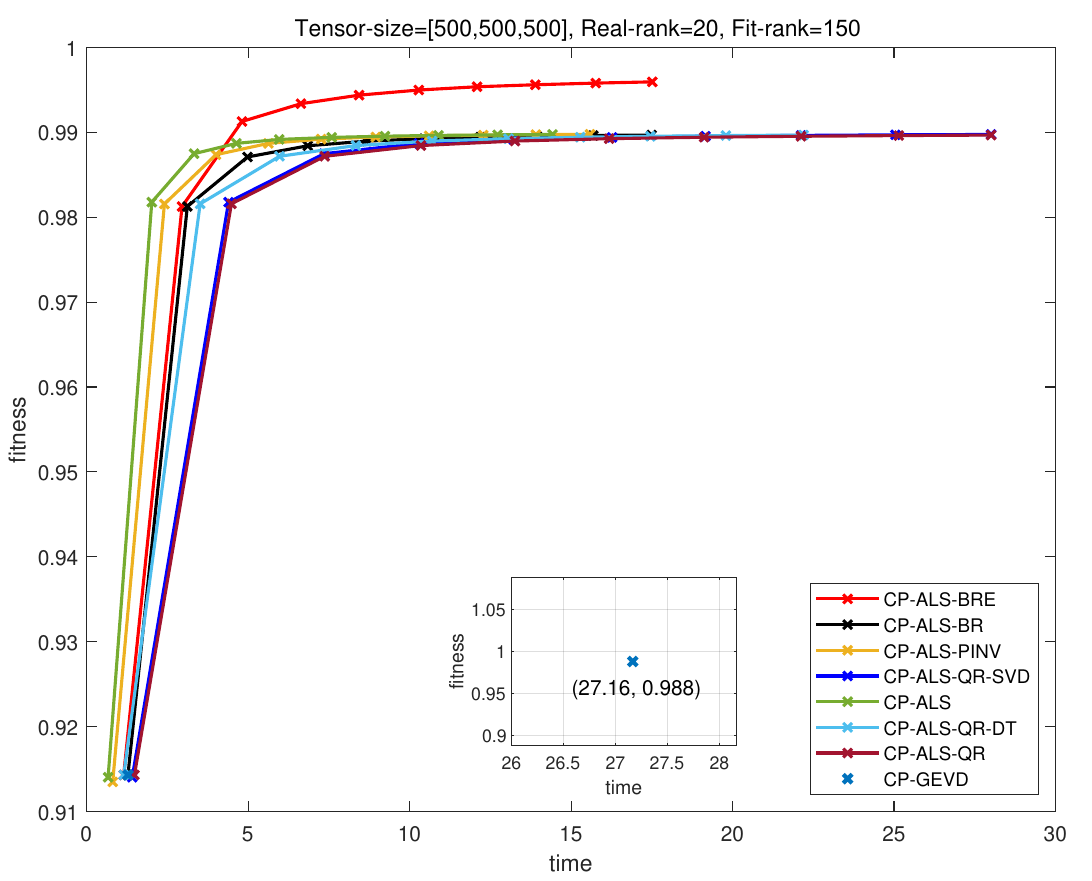}}
\end{minipage}
\begin{minipage}{0.5\textwidth}
\centering
\subfloat[]{
\includegraphics[scale=0.38]{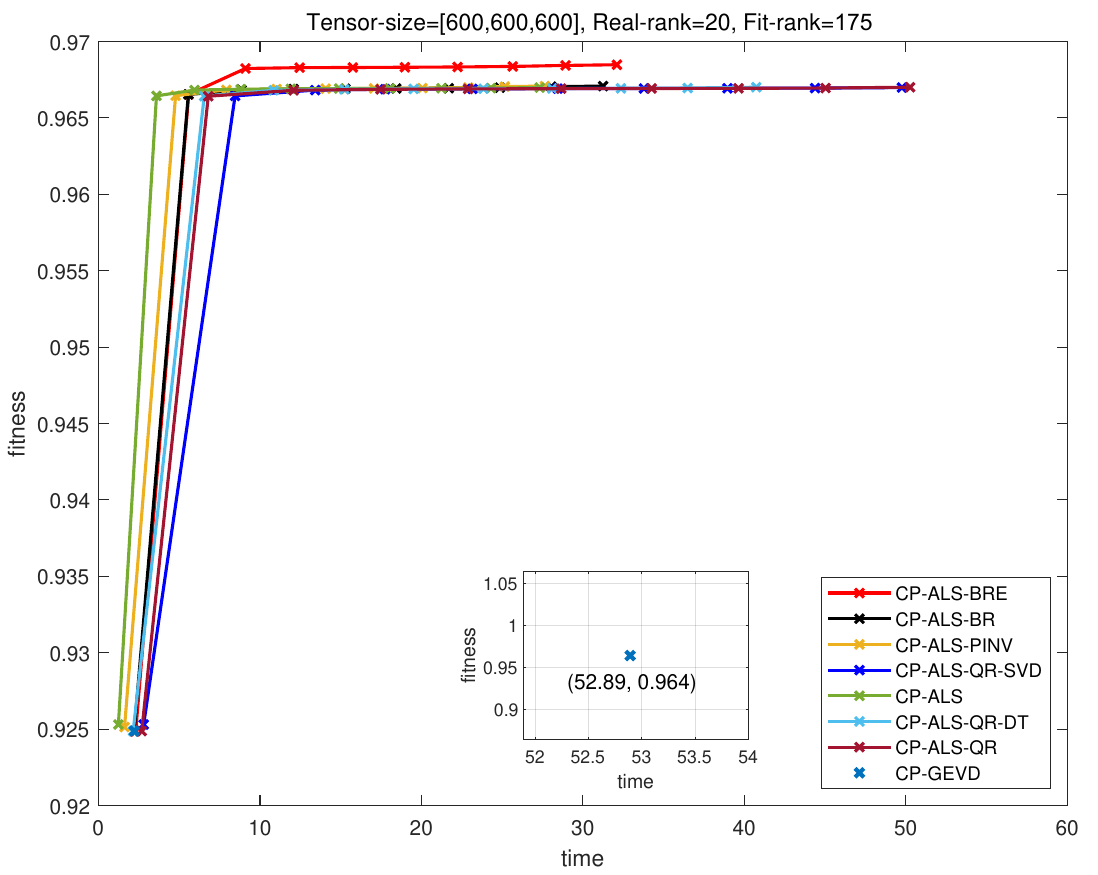}}
\end{minipage}
\caption{The performance plots of the third-order tensor. The left panel displays the performance for different fitting ranks and a single type of noise on the tensor of size [500, 500, 500]. The right panel shows the performance for different collinearity coefficients, different fitting ranks and mixed noise types on the tensor of size [600, 600, 600]. }
\label{fig-8}
\end{figure}

\begin{figure}[!ht]
\begin{minipage}[p]{0.5\textwidth}
\centering  
\subfloat[]{
\label{Fig9.sub.1}
\includegraphics[scale=0.38]{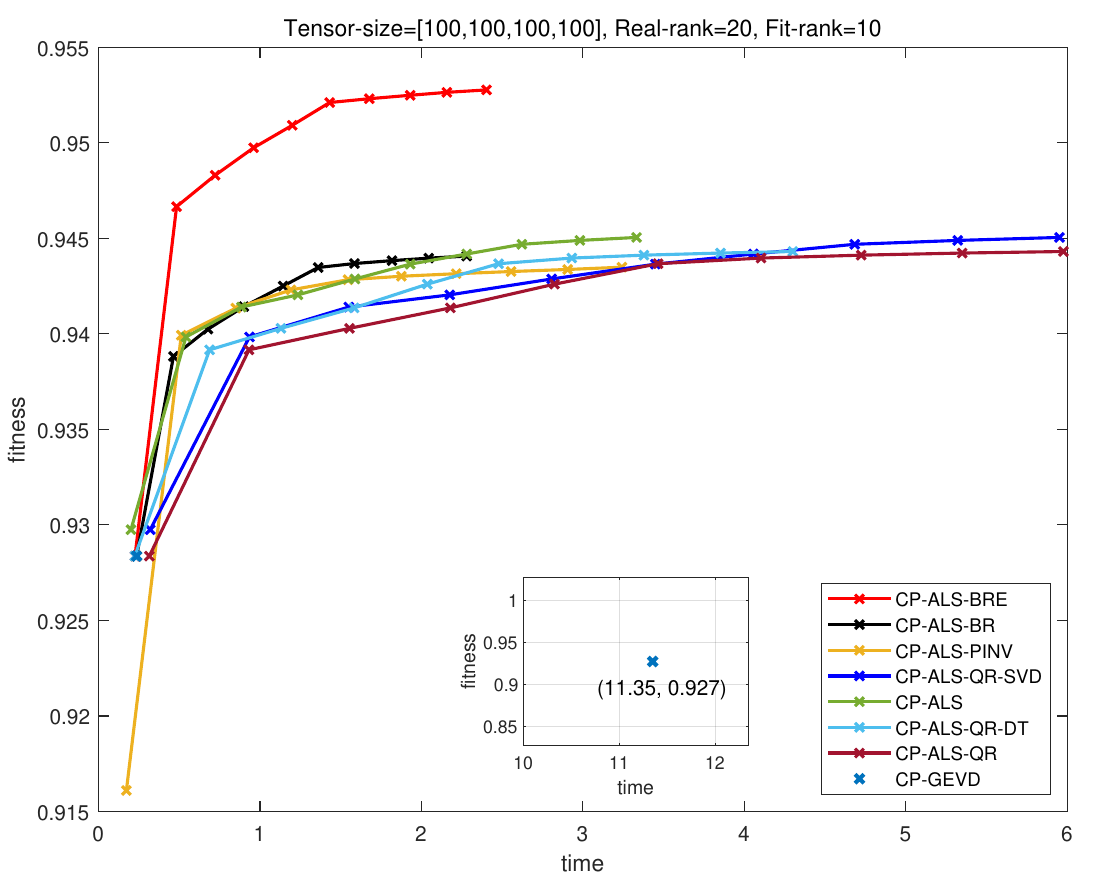}}
\vspace{10pt}
\end{minipage}
\begin{minipage}[p]{0.5\textwidth}
\centering
\subfloat[]{
\label{Fig9.sub.2}
\includegraphics[scale=0.38]{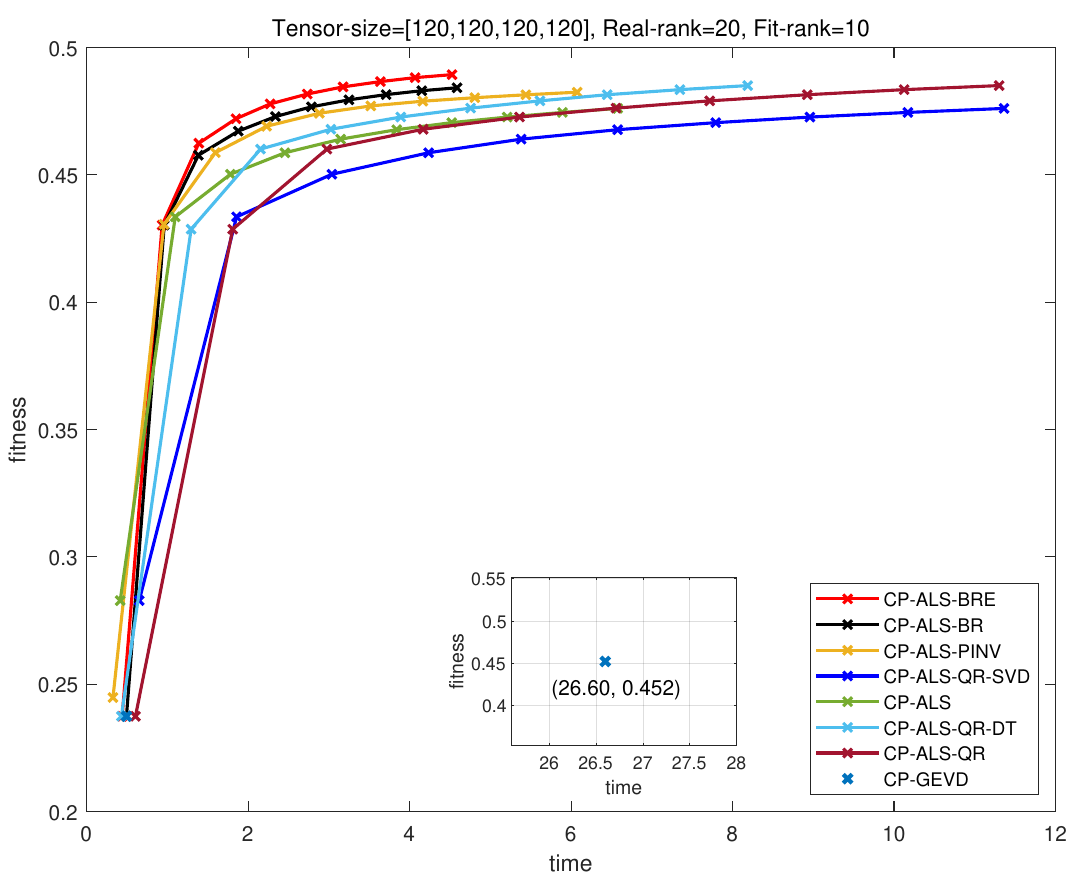}}
\vspace{10pt}
\end{minipage}
\begin{minipage}{0.5\textwidth}
\centering
\subfloat[]{
\label{Fig9.sub.3}
\includegraphics[scale=0.38]{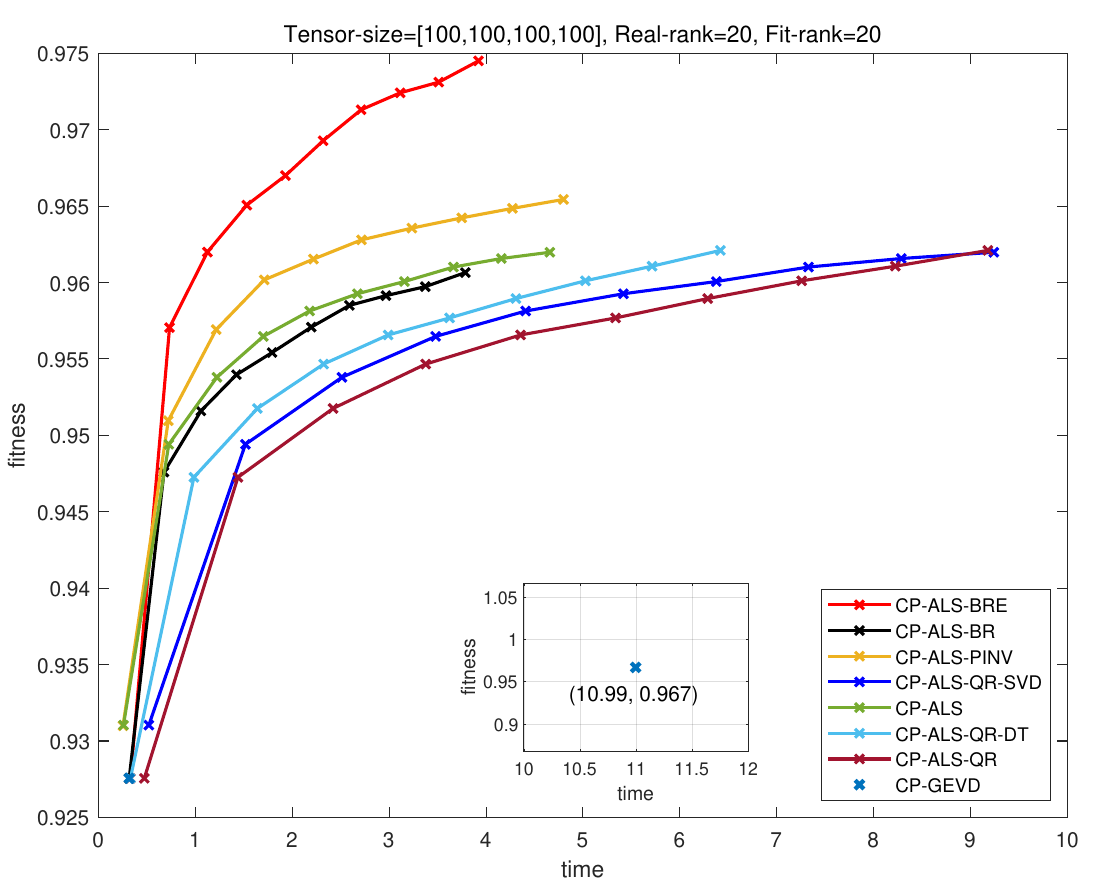}}
\vspace{10pt}
\end{minipage}
\begin{minipage}{0.5\textwidth}
\centering
\subfloat[]{
\label{Fig9.sub.4}
\includegraphics[scale=0.38]{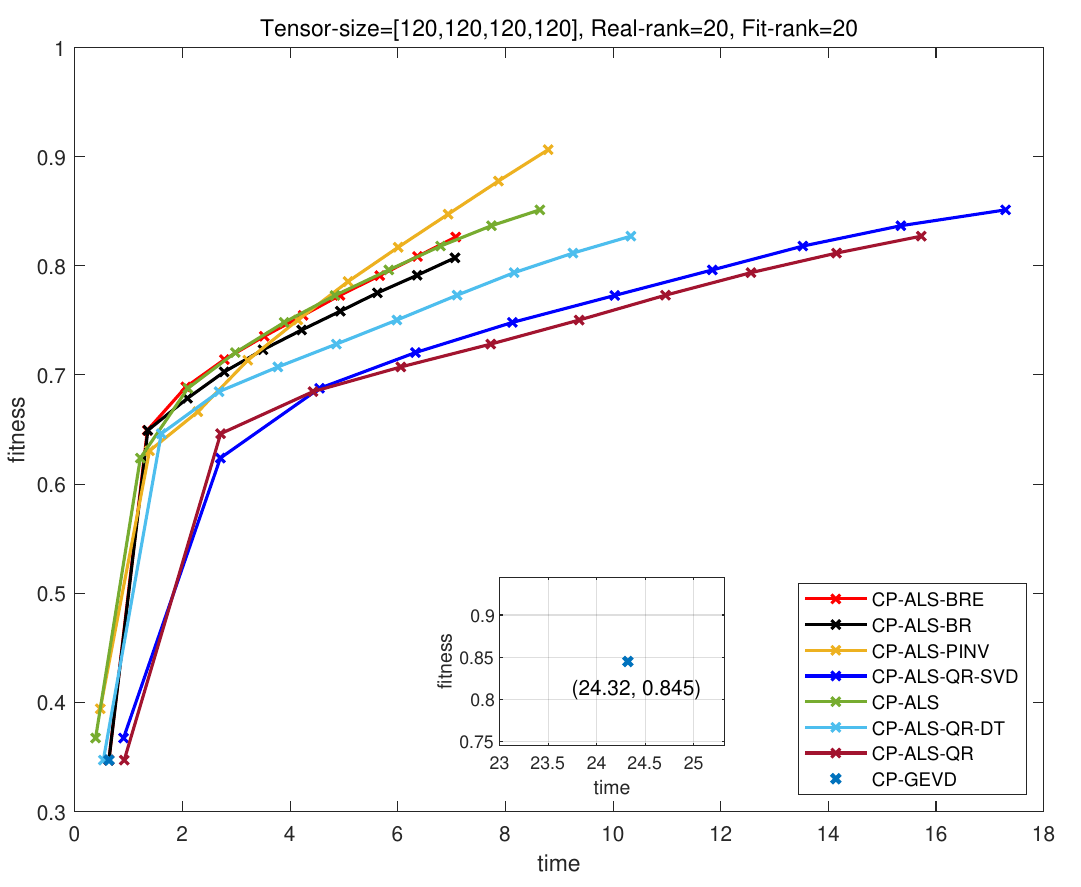}}
\vspace{10pt}
\end{minipage}
\begin{minipage}{0.5\textwidth}
\centering
\subfloat[]{
\label{Fig9.sub.5}
\includegraphics[scale=0.38]{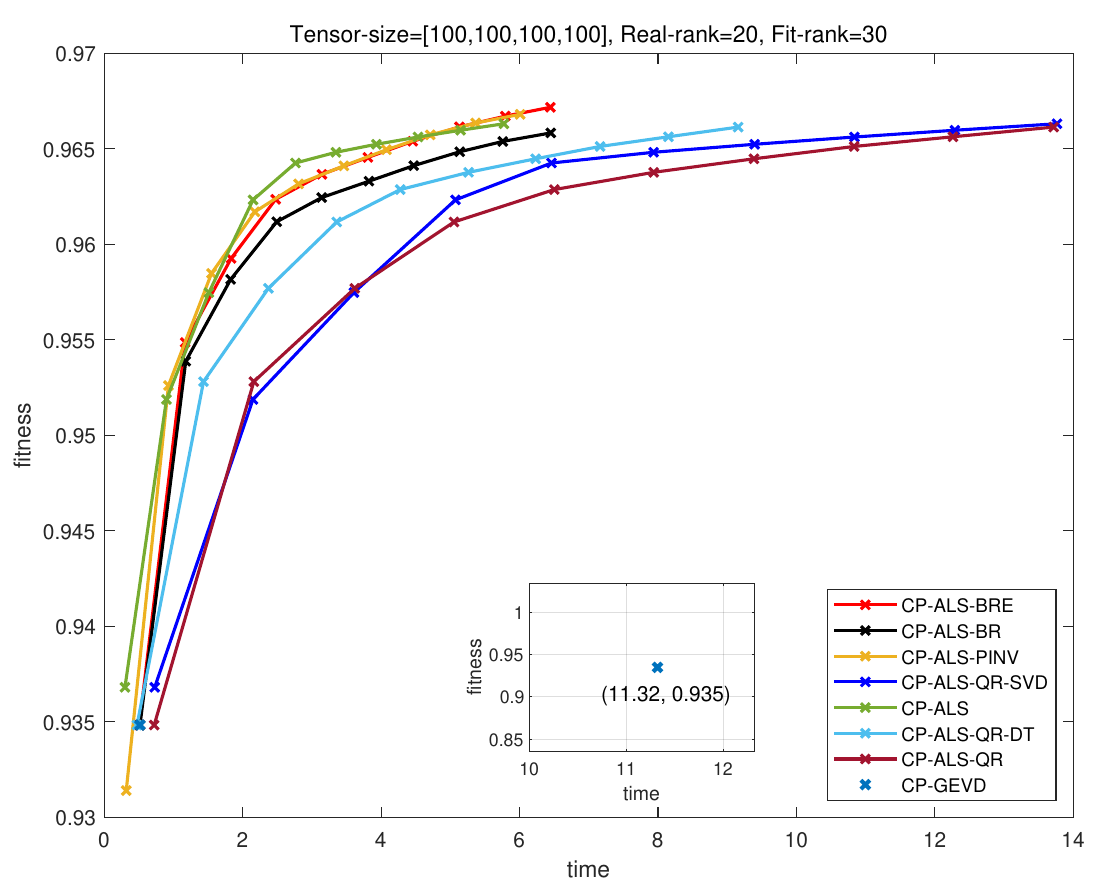}}
\end{minipage}
\begin{minipage}{0.5\textwidth}
\centering
\subfloat[]{
\label{Fig9.sub.6}
\includegraphics[scale=0.38]{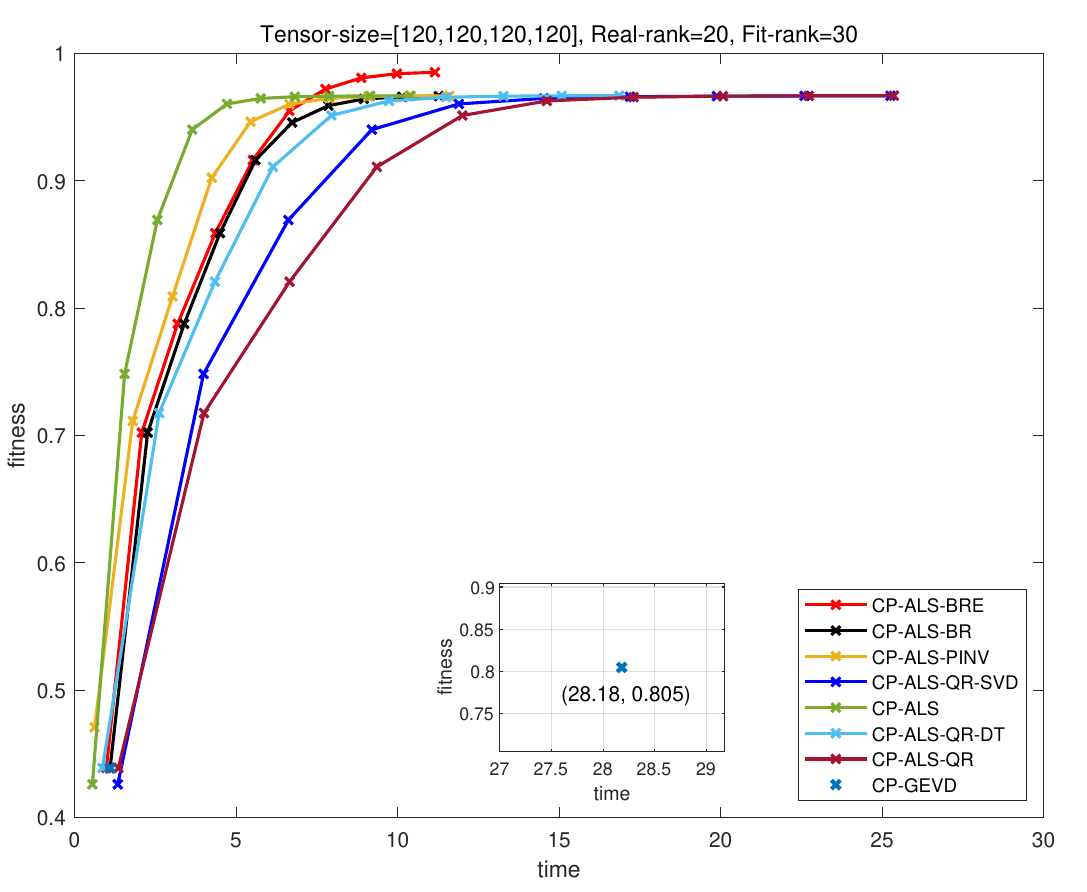}}
\end{minipage}
\caption{The performance plots of the fourth-order tensor. The left panel displays the performance for different fitting ranks and a single type of noise on a tensor of size [100, 100, 100, 100]. The right panel shows the performance for different collinearity coefficients, different fitting ranks, and mixed noise types on a tensor of size [120, 120, 120, 120].}
\label{fig-9}
\end{figure}

We present the experimental results of ALS-QR-BRE, ALS-QR-BR,CP-ALS-PINV, CP-ALS-QR-SVD, CP-ALS, CP-ALS-QR-DT, CP-ALS-QR and CP-GEVD  algorithms on synthetic third-order and fourth-order tensors in Figures \ref{fig-8} and \ref{fig-9}. The experiments involve running each algorithm (excluding the CP-GEVD algorithm) for 20 iterations, after which we compare the time required for 20 iterations and the fitting accuracy. Since CP-GEVD is a non-iterative algorithm, we compare its runtime and fitting accuracy with those of iterative algorithms after 20 iterations. From Figure \ref{fig-8}, we conclude that our proposed ALS-QR-BRE and ALS-QR-BR algorithms require approximately half the iteration time of the CP-ALS-QR and CP-ALS-QR-SVD algorithms across different fitting ranks. Although the CP-ALS-QR-DT algorithm accelerates the CP-ALS-QR algorithm, the improvement is minimal, and our proposed algorithms remain faster than CP-ALS-QR-DT. In most cases, our proposed algorithms exhibit comparable or slightly faster iteration speeds than CP-ALS and CP-ALS-PINV, as shown in Figures \ref{fig-8}(b) and \ref{fig-8}(d). Figures \ref{fig-8}(a)-(e) show that the ALS-QR-BRE algorithm has a significant runtime advantage over the CP-GEVD algorithm. Moreover, the ALS-QR-BRE algorithm consistently achieves superior fitting accuracy compared to the other algorithms across all synthesized third-order tensors and fitting ranks. From Figure \ref{fig-9}, we find that our proposed ALS-QR-BRE and ALS-QR-BR algorithms require less than half the iteration time of CP-ALS-QR and CP-ALS-QR-SVD across different fitting ranks. While the CP-ALS-QR-DT algorithm accelerates CP-ALS-QR and CP-ALS-QR-SVD, it still lags behind our proposed algorithms. For the majority of cases, oour proposed algorithm matches or outperforms CP-ALS and CP-ALS-PINV in iteration speed, except for a slight delay in Figures \ref{fig-9}(e) and \ref{fig-9}(f). Across different fitting ranks, ALS-QR-BRE consistently completes 20 iterations faster than the total runtime of CP-GEVD and generally achieves higher approximation accuracy.

The synthetic tensor experiments (Figures \ref{fig-8} and \ref{fig-9}) demonstrate that the proposed  ALS-QR-BRE algorithm significantly reduces the iteration time of CP-ALS-QR and outperforms CP-ALS-QR-DT in computational efficiency. Its iteration speed is generally comparable to or faster than CP-ALS, significantly faster than CP-GEVD, while its fitting accuracy outperforms other algorithms in most cases.

\subsection{Real-world datasets}\label{real_tensor}
In this subsection, we evaluate our algorithm on five real-world datasets: \emph{Indian Pines}, \emph{Salines}, \emph{Density}, \emph{Tabby Cat}, and \emph{Winter Landscape}. 

\begin{itemize}
\item 
\emph{Indian Pines}\footnote{\url{https://www.ehu.eus/ccwintco/index.php/Hyperspectral_Remote_Sensing_Scenes}}: The dataset is hyperspectral data collected by the AVIRIS sensor from the Indian Pines test site in northwestern Indiana. The data contains $145\times145$ pixels and 224 spectral bands, primarily covering agricultural land and natural vegetation such as forests.

\item \emph{Salines}: This dataset is a hyperspectral remote sensing image also collected by the AVIRIS over the Salinas Valley region in California, USA. The dataset contains 224 contiguous spectral bands with high spatial resolution. The spatial dimensions span 512 rows and 217 columns, forming a third-order tensor with dimensions $512\times217\times224$.

\item 
\emph{Density}\footnote{\url{https://gitlab.com/tensors/tensor_data_miranda_sim}}: 
This dataset is a third-order tensor representing the dynamic interaction between two fluids with different densities in a simulation. It has a 3D grid of size $2048 \times 256 \times 256$ (Z, Y, X), stored in double precision format (1 GB), providing detailed insights into fluid dynamics.
\item
\emph{Tabby Cat}\footnote{\url{https://www.pexels.com/video/video-of-a-tabby-cat-854982/}}: This dataset is derived from a standard color video source, where the initial temporal segment of the footage has been processed to form a fourth-order tensor with dimensions $720 \times 1280 \times 3 \times 100$. 

\item 
\emph{Winter Landscape}\footnote{\url{https://www.pexels.com/video/an-animal-on-a-winter-landscape-2313069/}}: This dataset also is derived from a standard color video source, where the initial temporal segment of the footage has been processed to construct a fourth-order tensor with dimensions $1510 \times 1080 \times 3\times 60$. 
\end{itemize}

\begin{figure}[!ht]
 \begin{minipage}[p]{0.5\textwidth}
\centering  
\subfloat{
\includegraphics[scale=0.38]{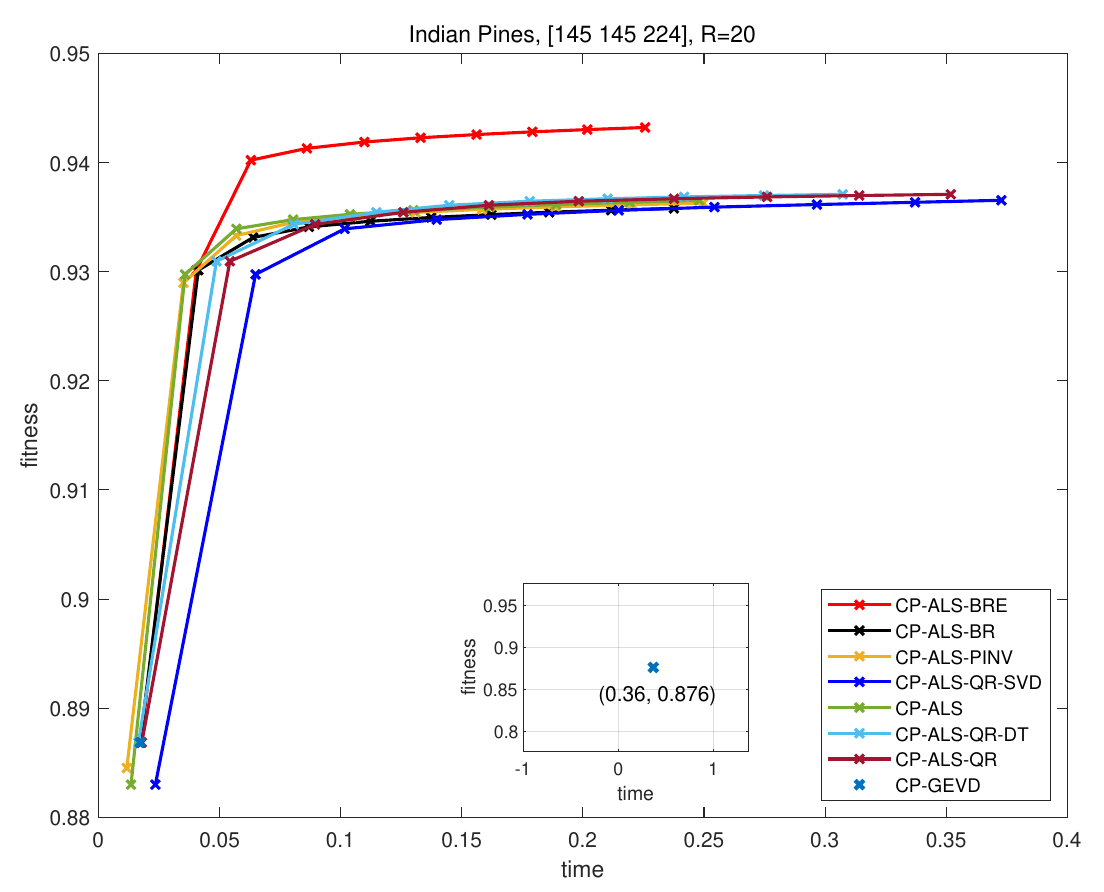}}
\end{minipage}
\begin{minipage}[p]{0.5\textwidth}
\centering
\subfloat{
\includegraphics[scale=0.373]{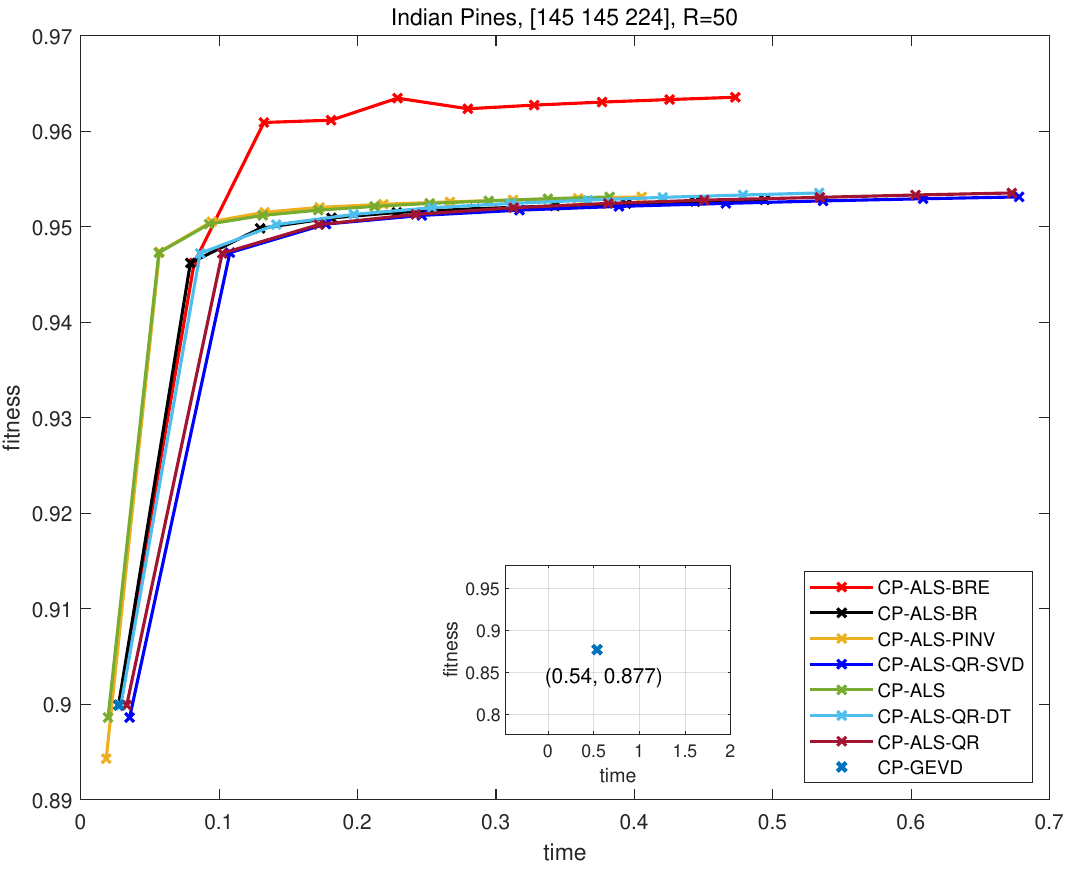}}
\end{minipage}
\vspace{-5pt}
\caption{Performance graphs of each algorithm on the Indian Pines dataset. The left and right panels display the performance plots for different algorithms with ranks of 20 and 50, respectively.}
\label{fig_10}
\end{figure}

\begin{figure}[!ht]
\vspace{-20pt}
\begin{minipage}[p]{0.5\textwidth}
\centering  
\subfloat{
\includegraphics[scale=0.38]{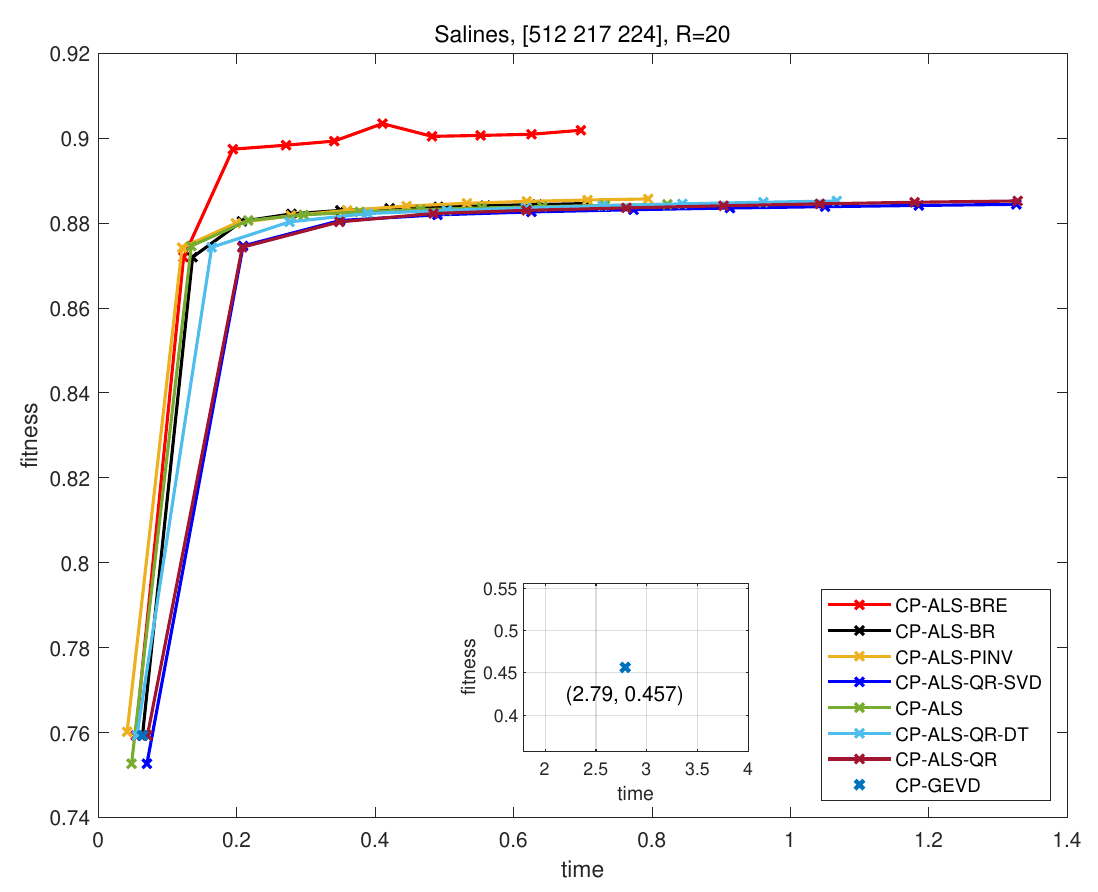}}
\end{minipage}
\begin{minipage}[p]{0.5\textwidth}
\centering
\subfloat{
\includegraphics[scale=0.38]{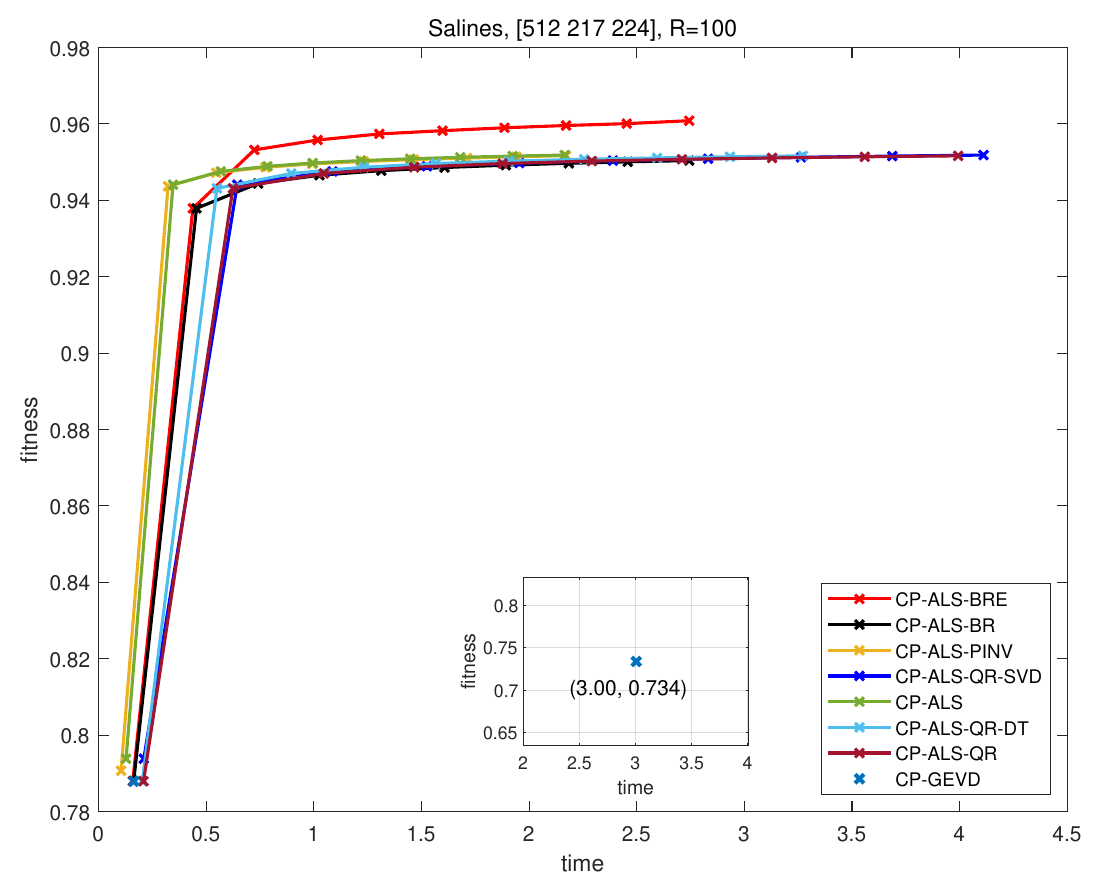}}
\end{minipage}
\vspace{-5pt}
\caption{Performance plots of each algorithm on the Salines dataset. The left panel displays the performance plot with a fitting rank of 20, while the right panel shows the performance plot with a fitting rank of 100.}
\label{fig_11}
\end{figure}

\begin{figure}[!ht]
 \begin{minipage}[p]{0.5\textwidth}
\centering  
\subfloat{
\includegraphics[scale=0.38]{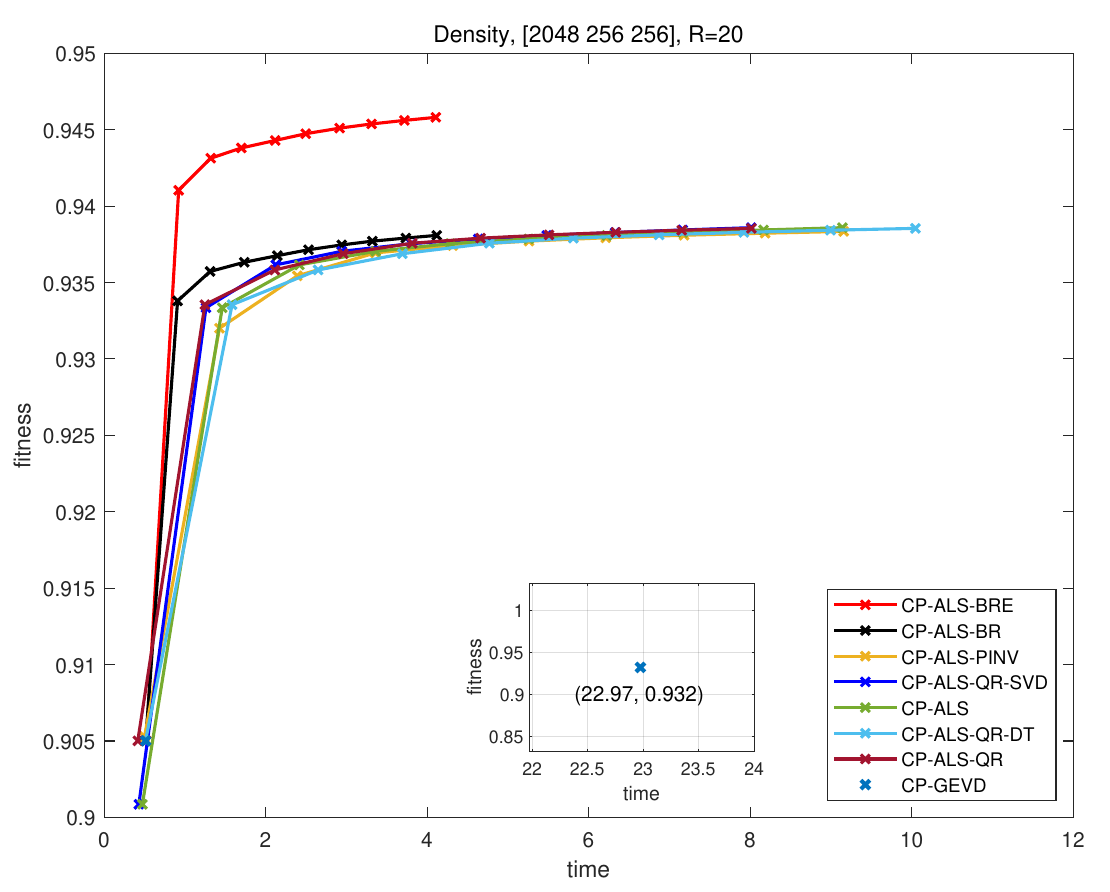}}
\end{minipage}
\begin{minipage}[p]{0.5\textwidth}
\centering
\subfloat{
\includegraphics[scale=0.38]{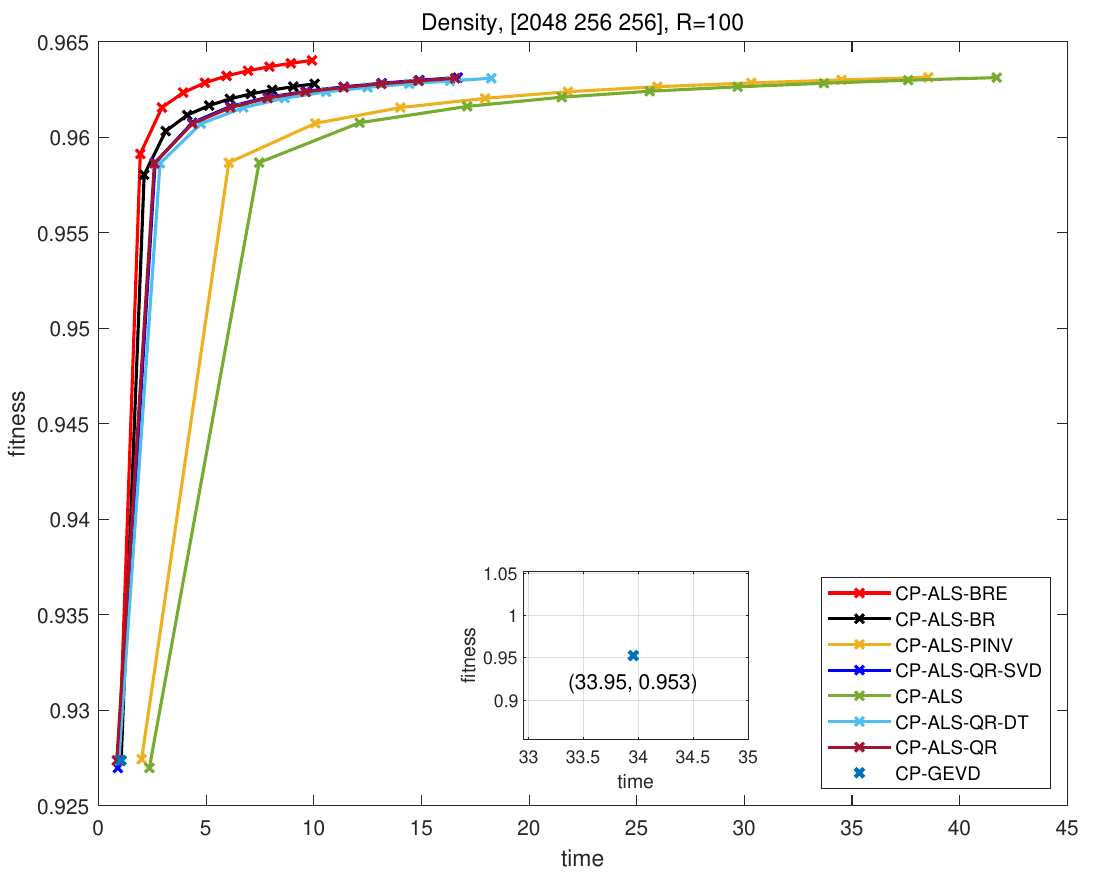}}
\end{minipage}
\vspace{-5pt}
\caption{Performance plots of each algorithm on the Density dataset. The left panel displays the performance plot with a fitting rank of 20, while the right panel shows the performance plot with a fitting rank of 100.}
\label{fig_12}
\end{figure}

\begin{figure}[!ht]
\begin{minipage}{0.5\textwidth}
\centering
\subfloat{
\includegraphics[scale=0.38]{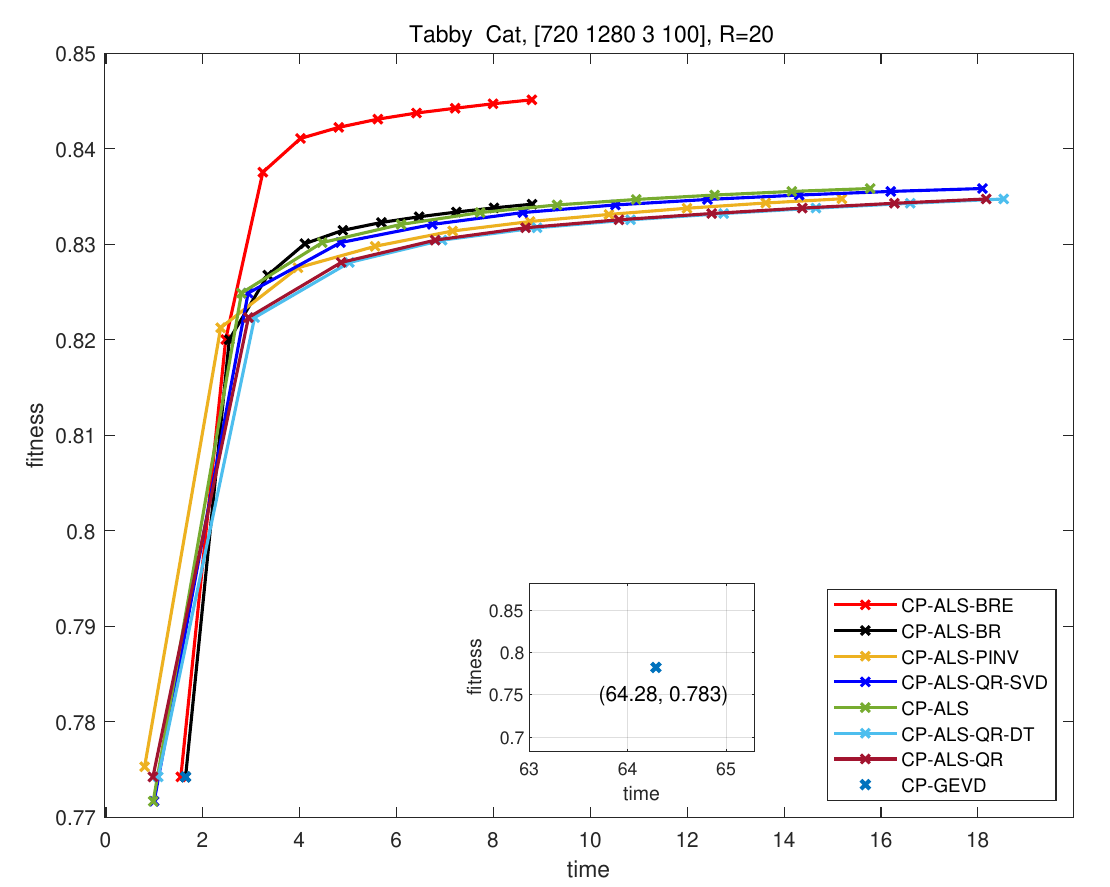}}
\end{minipage}
\begin{minipage}{0.5\textwidth}
\centering
\subfloat{
\includegraphics[scale=0.38]{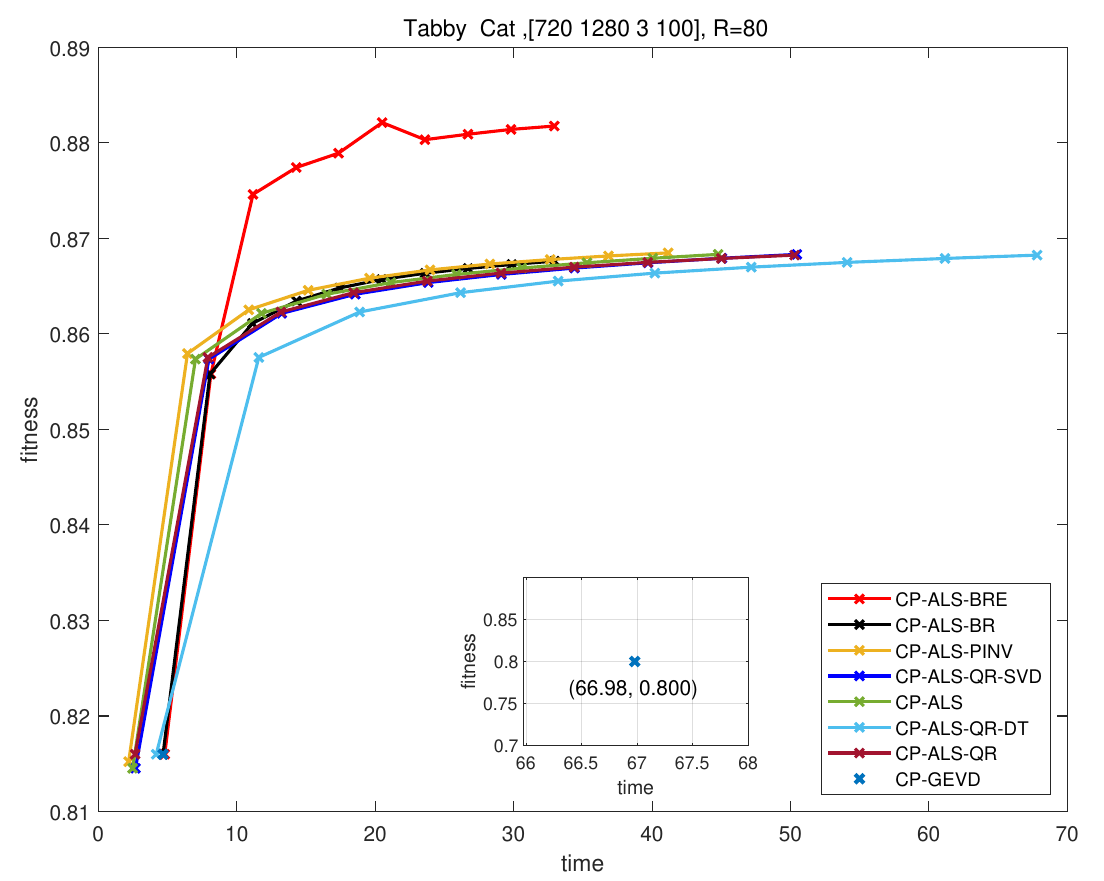}}
\end{minipage}
\caption{Plots illustrating the performance of each algorithm on the Tabby Cat dataset. The left and right panels display the performance plots for different algorithms with ranks of 20 and 80, respectively. }
\label{fig_13}
\end{figure}

\begin{figure}[!ht]
\begin{minipage}{0.5\textwidth}
\centering
\subfloat{
\includegraphics[scale=0.38]{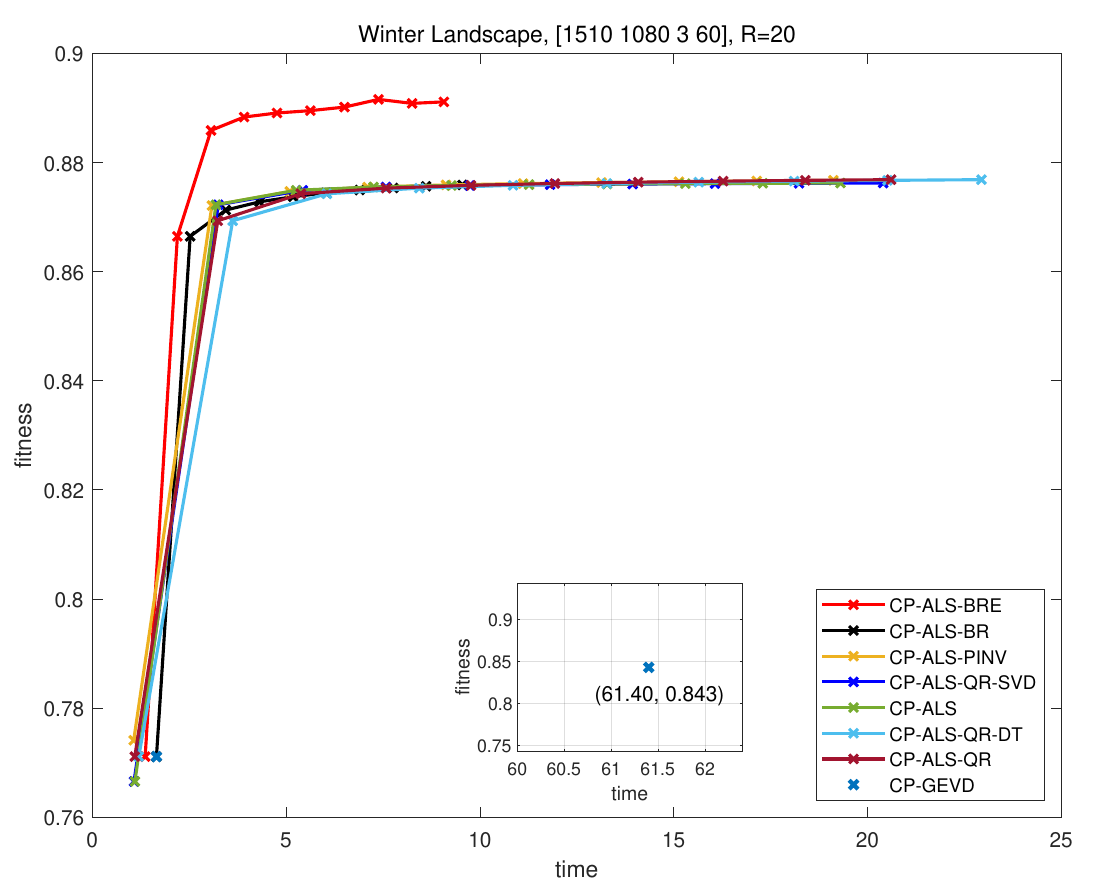}}
\end{minipage}
\begin{minipage}{0.5\textwidth}
\centering
\subfloat{
\includegraphics[scale=0.38]{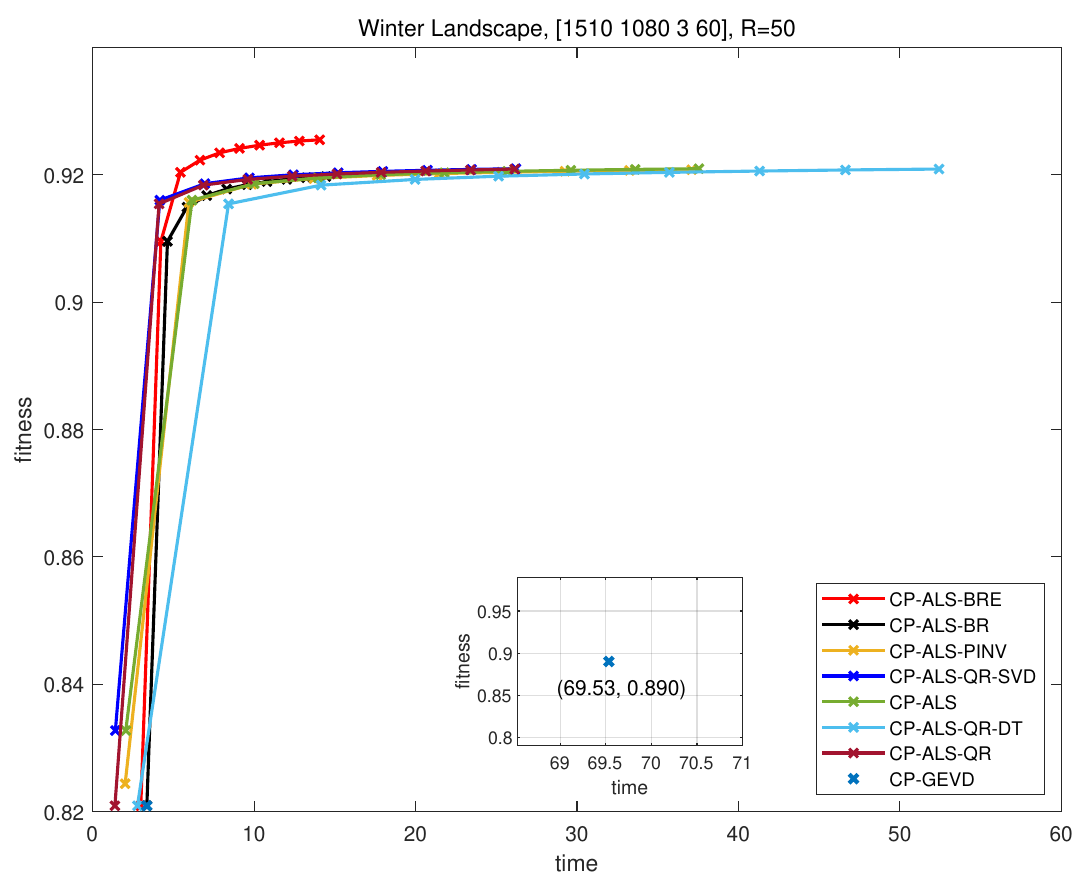}}
\end{minipage}
\caption{Plots illustrating the performance of each algorithm on Winter Landscape dataset. The left and right panels display the performance plots for different algorithms with ranks of 20 and 50, respectively. } 
\label{fig_14}
\end{figure}

For each real-world tensor experiment, we run the ALS-QR-BRE, ALS-QR-BR, CP-ALS-PINV, CP-ALS-QR-SVD, CP-ALS, CP-ALS-QR-DT, CP-ALS-QR and CP-GEVD algorithms to evaluate and compare their operational efficiency and fitting accuracy. In the experiment, all iterative algorithms are run for 20 iterations, with results presented in Figures \ref{fig_10}-\ref{fig_14}. Since the computational overhead introduced by extrapolation is negligible, the time advantage of the ALS-QR-BRE iterations is essentially the same as that of the ALS-QR-BR algorithm.

We present the rank-20 and rank-50 decompositions of the third-order tensor datasets \emph{Indian Pines} in Figure \ref{fig_10}. It shows that for the rank-20 decomposition on the \emph{Indian Pines} dataset, Our proposed ALS-QR-BRE algorithm offers faster iteration speed and improves fitting accuracy by 0.79\% compared to other iterative methods, significantly outperforming the CP-GEVD algorithm. For rank-50 decomposition, it maintains computational advantages over all algorithms except CP-ALS and CP-ALS-PINV, improving fitting accuracy by 1.13\% over other iterative methods and achieving much higher accuracy than CP-GEVD.  

For the \emph{Salines} dataset, we conduct rank-20 and rank-100 decompositions in Figure \ref{fig_11}. At rank 20, ALS-QR-BRE achieves the highest efficiency and 1.87\% better fitting accuracy than other iterative methods, significantly surpassing CP-GEVD. At rank 100, it maintains computational advantages over all but CP-ALS and CP-ALS-PINV, with a 1.06\% accuracy improvement over iterative methods and notably better performance than CP-GEVD.

We perform rank-20 and rank-100 decompositions for the \emph{Density} dataset, and the results are shown in Figure \ref{fig_12}. The results show that for rank-20 decomposition, our proposed ALS-QR-BRE algorithm completes iterations in 52.25\% of the time of CP-ALS-QR and CP-ALS-QR-SVD, 45.83\% of the time of CP-ALS and CP-ALS-PINV, and 19.19\% of the runtime of CP-GEVD, while improving fitting accuracy by at least 0.79\% compared to other algorithms. For rank-100 decomposition, ALS-QR-BRE requires only 60.33\% of the time of CP-ALS-QR and CP-ALS-QR-SVD, 24.21\% of CP-ALS and CP-ALS-PINV, and 31.20\% of CP-GEVD, achieving higher fitting accuracy than the other methods. 

Finally, in the experiment on real-world datasets, we use two fourth-order tensor datasets, \emph{Tabby Cat} and \emph{Winter Landscape}, to demonstrate the performance of our algorithm on large-scale tensors. We execute rank-20 and rank-80 decompositions on the \emph{Tabby Cat} dataset (results in Figure \ref{fig_13}) and rank-20 and rank-50 decompositions on the \emph{Winter Landscape} dataset (results in Figure \ref{fig_14}). From Figure \ref{fig_13}, for rank-20 decomposition, our proposed ALS-QR-BRE algorithm improves fitting accuracy by at least 1.11\% compared to others and completes iterations in 49.06\% of the time of CP-ALS-QR, 56.62\% of CP-ALS, and 14.62\% of CP-GEVD. For rank-80 decomposition, it finishes in 65.35\% of CP-ALS-QR time, 73.57\% of CP-ALS time, and 51.75\% of CP-GEVD time, with at least 1.53\% accuracy improvement.  Figure \ref{fig_14} shows that for rank-20 decomposition, the proposed ALS-QR-BRE enhances fitting accuracy by 1.63\% and requires 44.93\%, 48.09\%, and 15.87\% of the iteration time of CP-ALS-QR, CP-ALS, and CP-GEVD respectively. For rank-50 decomposition, it maintains this speed advantage, completing iterations in 53.92\% of CP-ALS-QR time, 37.83\% of CP-ALS time, and 21.50\% of CP-GEVD time, while achieving higher fitting accuracy. 

Experimental evaluations on synthetic and real-world tensor datasets show that the proposed ALS-QR-BRE algorithm consistently achieves shorter iteration times than CP-ALS-QR and CP-ALS-QR-SVD. When the first one or two dimensions of the real-world tensor data are significantly larger than the subsequent dimensions, our algorithm’s iteration speed is significantly faster than CP-ALS, CP-ALS-PINV, and CP-GEVD. Furthermore, in synthetic tensor applications, ALS-QR-BRE outperforms other compared algorithms in fitting accuracy in most cases and similarly demonstrates superior fitting performance on real-world tensor datasets.

\section{Conclusion}\label{conclusion}
In this paper, we introduced restructured dimension tree to improve the utilization of intermediate tensors, thereby reducing the number of TTM operations and computational complexity. Additionally, we designed customized extrapolation acceleration technique for the matrix $\mathbf{Q}_0$ in the CP-ALS-QR algorithm and presented its theoretical rationale. By integrating these techniques, we proposed the ALS-QR-BRE algorithm. Experiments on both synthetic and real-world tensor datasets demonstrated that ALS-QR-BRE outperformed existing algorithms in computational efficiency and fitting accuracy.

Future research can explore several potential optimization directions. When dealing with large-scale datasets and high-rank decompositions, the computation of $\mathbf{Q}_0$ remains the primary bottleneck. Exploiting its upper triangular structure may further reduce computational costs, which is particularly significant for large-scale datasets. Additionally, further research is needed to develop simple and effective acceleration methods under exact precision without increasing the computational load.

\section*{Declarations}
{\bf Funding:} This research is supported by the National Natural Science Foundation of China (NSFC) grants 92473208, 12401415, the Key Program of National Natural Science of China 12331011, the 111 Project (No. D23017), the Innovative Research Group Project of Natural Science Foundation of Hunan Province of China (No. 2024JJ1008), the Natural Science Foundation of Hunan Province (No. 2025JJ60009).

\noindent{\bf Data Availability:} Enquiries about data/code availability should be directed to the authors.

\noindent{\bf Competing interests:} The authors have no competing interests to declare that are relevant to the content of this paper.


\bibliographystyle{abbrv}
\bibliography{ALS_QR_BRE}


\end{document}